\documentstyle[12pt]{article}
\date{}

\newcommand{\be}{\begin{equation}}
      \newcommand{\ee}{\end{equation}}
      \newcommand{\ba}{\begin{eqnarray}}
       \newcommand{\ea}{\end{eqnarray}}
\newcommand{\ban}{\begin{eqnarray*}}
       \newcommand{\ean}{\end{eqnarray*}}

\newcommand{\lexp}{exp^{\cal{L}, \tau}_{p}}

 \newcommand{\qed}{\hspace*{\fill}\rule{3mm}{3mm}\quad}
 \newcommand{\Pf}{\noindent {\em Proof.} }

\newcommand{\sect}[1]{\section{#1} \setcounter{equation}{0}}

\newtheorem{lem}{Lemma}[section]

\begin{document}
\newtheorem{defn}[lem]{Definition}
 \newtheorem{theo}[lem]{Theorem}
 \newtheorem{prop}[lem]{Proposition}
 \newtheorem{rk}[lem]{Remark}
 \newtheorem{ex}[lem]{Example}
 \newtheorem{note}[lem]{Note}
 \newtheorem{conj}[lem]{Conjecture}
 \newtheorem{cor}[lem]{Corollary}

\title{On the $l$-Function and the Reduced Volume of Perelman I
}
\author{Rugang Ye\\Department of Mathematics\\University of California,
Santa Barbara} 
\maketitle

%\noindent Table of Contents \\
%Section 1 Introduction \,  $\cdot \cdot \cdot \cdot$ Page 1
%\\
%Section 2 \, Basic Properties of the $l$-Function
%I  $\cdot \cdot \cdot \cdot$ Page 3 \\
%Section 3 \, Basic Properties of the $l$-Function II  $\cdot \cdot \cdot \cdot$ Page 17 \\
%Section 4 \, The Reduced Volume  $\cdot \cdot \cdot \cdot$ Page 20 \\
%Section 5 \, Asymptotic Limits of $\kappa$-Solutions $\cdot \cdot \cdot \cdot$ Page 27 \\
%\\
%\vspace{1cm}

\sect{Introduction}

In [P1], Perelman introduced, among other things, two important
tools for analyzing the Ricci flow: the reduced distance, i.e. the
$l$-function, and the reduced volume.  The $l$-function is defined
in terms of a natural curve energy along the Ricci flow, which is
analogous to the classical curve energy employed in the study of geodesics,
but involves the evolving metric, as well as the scalar curvature as a
potential term. The reduced volume is a certain integral
involving the $l$-function. The $l$-function and the reduced
volume enjoy a number of very nice analytic and geometric properties, 
including in particular the fundamental  monotonicity of the reduced 
volume.
 These properties can be used, as demonstrated by Perelman,
to classify and analyze blow-up limits of the Ricci flow, and to
obtain various estimates for the Ricci flow, such as
non-collapsing estimates and curvature estimates. 
%A notable
%example is Prop. 11.2 in [P1], which identifies blow-down limits of
%$\kappa$-solutions (non-collapsed ancient solutions of the
%Ricci flow with bounded curvature and nonnegative curvature
%operator) to be nonflat gradient shrinking solitons. This result plays an important role in Perelman's anlysis of blow-up 
%singularities of the Ricci flow in [P1] and [P2]. 

The main purpose of this paper is to  present a number of analytic
and geometric properties of the $l$-function and the reduced
volume, including in particular the monotonicity, the upper bound and the
rigidities of the reduced volume. In Perelman's paper, a general
assumption concerning the $l$-function and the reduced volume is
uniformly bounded sectional curvature. The results obtained in [P1]
under this assumption are sufficient for the application to the
geometrization of 3-manifolds in [P2]. Because of the fundamental role of
the $l$-function and the reduced volume for analyzing Ricci flow
in general,  it is very desirable to allow weaker geometric
conditions. Our main focus is to deal with the situation in which
only a lower bound for the Ricci curvature is assumed. On the
other hand, we hope that our treatment can provide assistance for
understanding Perelman's theory, even when one is only interested
in the  case of bounded sectional curvature.

For the convenience of the reader, we give here a short account of the 
main topics in this paper. In Section 2, we first present the basic concepts such as the 
$l$-function, the $\cal L$-geodesics, and  the $\cal L$-exponential map.
Their basic properties are then analysed, which include in particular  
the local Lipschitz properties   and  
the local semi-concavity of the $l$-function. Several basic estimates 
for the $l$-function are also presented.  One highlight of this section is 
Theorem 2.20, which establishes the fundamental differential inequalities 
(2.67) and (2.68) everywhere in the weak sense, which are formulated in 
(2.69) and (2.70). Another highlight is the Harnack inequality for the 
$l$-function in Theorem \ref{3estimateslemma}. 

In Section 3, we first present an estimate for the minimum value of the $l$-function.
Next we present an important estimate of Perelman  which 
provides a lower bound  for the $l$-function at any time in terms of the squared distance at the 
same time
from a fixed point, assuming nonnegative curvature operator, see Lemma 3.2. (A corresponding 
upper bound is also included.)  
In comparison, the easier Lemma 2.3 provides a similar estimate in terms of 
the squared distance at an earlier time (assuming a lower bound for the Ricci curvature)
or a later time (assuming an upper bound for the Ricci curvature).  Perelman's estimate provides 
an important analytic ingredient for dealing with some crucial and delicate convergence issues
of integrals involving the $l$-function. 
%The most notable application of this estimate is 
%perhaps for proving [Proposition 11.2, P1], see [Y2] (or [Y1]).
%This estimate, along with its counterpart Lemma 2.2, is also used for proving  
Another topic in this section is Theorem 3.3, which expands the scope of  Theorem 2.20 to admit test functions 
which may not have compact support but satisfy a certain decay condition. 
Such expansions are needed for applying the differential inequalities 
(2.67) and (2.68) to topics around the reduced volume, as in Section 4 and 
[Y2] (see also [Y1]).  

In Section 4, we derive a number of important properties of the reduced 
volume, which include the monotonicity (Theorem 4.5), the upper bound (Theorem 4.3), and the rigidity
regarding the upper bound
(Theorem 4.4). We also derive the 
rigidity regarding the monotonicity (Theorem \ref{solitontheorem}), which 
says that a solution of the backward Ricci flow  must be a gradient shrinking soliton if 
the values of the reduced volume are equal at two different times. This theme 
appears again in [Y2] (see also [Y1]) in a different set-up.   

%Section 5 focuses on Theorem 5.1, which is precisely  Proposition 11.2 in [P1]. 
%The proof of Theorem 5.1 follows basic ideas of Perelman, and is based on the analytic and
%geometric properties of the $l$-function and the reduced volume presented in the 
%previous sections. In particular, the monotonicy of the reduced volume, 
%the derivative estimates for the $l$-function, the soliton analysis in 
%Section 4, and the bounds for the $l$ function in Lemma 3.2 all 
%come into play in the proof. 

Communications with Perelman were of great help for understanding
his ideas.  We also
benefited much from conversations with Guofang Wei.  We would like to 
thank  
Vitali
Kapovich for helping to find the reference [GW].

This paper is part of [Y1], whose first version was posted on the author's webpage 
in February 2004.

\sect{Basic Properties of the $l$-Function I}

Consider a smooth solution $(M, g=g(\tau))$ of the backward Ricci flow
\be  \label{Ricciflow}
\frac{\partial g}{\partial \tau}=2Ric
\ee
on a manifold $M$ over an interval $[0, T)$.  We assume that $(M, g(\tau))$ is complete
for each $\tau \in [0, T)$.  Note that the theory presented here is meant to 
be applied to solutions of the Ricci flow. Indeed, a solution of the Ricci flow can be 
converted into a solution of the backward Ricci flow by a time reversal.  \\

\noindent {\bf Notations} We shall denote the distance between two
points $q_1, q_2$ with respect to the metric $g(\tau)$ by $d(q_1,
q_2, \tau)$, $d_g(q_1,q_2, \tau)$ or $d_{g(\tau)}(q_1, q_2)$. The geodesic ball of
center $q$ and radius $r$ with respect to the metric $g(\tau)$
will be denoted by $B_r(q, \tau)$. The volume form of $g(\tau)$
will be denoted by $dq$ or $dq|_{\tau}$. The scalar curvature
$R_{g(\tau)}$ of $g(\tau)$ at a point $q$ will be written as $R(q,
\tau)$. Similar notations
are also used for other curvature quantities.\\

A basic and simple lemma  is this.

\begin{lem} If $Ric\geq -cg$ for a nonnegative constant $c$ on the time
interval $[0, \tau]$, then
 \be \label{basicboundI} e^{-2cs}g(0)\le g(s) \le
e^{2c(\tau-s)} g(\tau) \ee
 for $s \in [0, \tau]$. If
$Ric \le Cg$ for a nonnegative constant $C$ on $[0, \tau]$, then
 \be \label{basicboundII} e^{2C(s-\tau)}g(\tau)\le g(s) \le
e^{2Cs}g(0) \ee for $s \in [0, \tau]$.
\end{lem}

We consider Perelman's $\cal L$-energy for piecewise $C^1$ curves
$\gamma: [a, b] \rightarrow M, 0\leq a <b <T$: \be \label{Lenergy}
{\cal{ L}}_{a, b}(\gamma)=\int_a^{b} \sqrt{s}(R(\gamma(s),
s)+|\dot{\gamma}|^2) ds, \ee where  $| \cdot |=|\cdot |_{g(s)}$.
For a given $\tau$ we abbreviate ${\mathcal{L}}_{0, \tau}$ to
$\mathcal{L}$.
%For a given $\gamma$, we set
%$X=\dot{\gamma}$.
The ${\cal{L}}_{a, b}$-geodesic (or $\cal{L}$-geodesic) equation
is: \be \label{Lgeodesic} \nabla_{\frac{d}{ds}}
\dot{\gamma}-\frac{1}{2}\nabla R+
\frac{1}{2s}\dot{\gamma}+2Ric(\dot{\gamma}, \cdot)=0,
 \ee
 where $R=R_{g(s)}, Ric=Ric_{g(s)}$, and $\nabla$ is the Levi-Civita
 connection of $g(s)$. This is the Euler-Lagrange equation of the
 $\cal L$-energy. Its (smooth) solutions are called ${\cal{L}}_{a,
 b}$-geodesics or $\cal{L}$-geodesics.

To better understand the properties of ${\cal{L}}_{a, b}$-geodesics, it is helpful to
introduce a convenient reparametrization. We set $t=\sqrt{s}$ and $\gamma'=d\gamma/dt
=2t\dot{\gamma}$.
Then
\be  \label{newL}
{\cal{L}}_{a, b}(\gamma)=\int_{\sqrt{a}}^{\sqrt{b}} (\frac{1}{2}|\gamma'|^2+2Rt^2)dt
\ee
and the ${\cal{L}}_{a, b}$-geodesic equation becomes
\be  \label{newgeodesic}
\nabla_{\frac{d}{dt}} \gamma'-2t^2\nabla R+4tRic(\gamma', \cdot)=0.
\ee

Next we choose a reference point $p\in M$ and define $L(q,
\tau)=L_g(q, \tau)$ to be the infimum of ${\cal L}(\gamma)$ for
$\gamma: [0, \tau] \rightarrow M$ with $\gamma(0)=p$ and $
\gamma(\tau)=q$. (We write $L_g(q, \tau)$ if we need to indicate
the dependence on the solution $g$.)\\

\noindent {\bf Definition 1}   We  define the {\it reduced
distance} (of Perelman) to be \be \label{lsmall} l(q, \tau)=l_g(q,
\tau)=\frac{L(q, \tau)}{2\sqrt{\tau}}.  \ee
We also call it 
the {\it $l$-function} (of Perelman).  The reference point 
$p$ will be called an {\it $l$-base}. \\

An easy computation leads to the following basic lemma.

\begin{lem} The $l$-function is invariant 
under the rescaling $g(\tau) \rightarrow g_a(\tau) \equiv a^{-1}g(a\tau)$, i.e.
\ba
l_{g_a}(q, \tau)=l_g(q, a\tau)
\ea
for all $\tau \in (0, \frac{T}{a})$ and $q \in M$.
\end{lem} 

Next we derive an estimate for $l$ in terms of the distance function.

\begin{lem}  \label{firstllemma} Assume that $Ric\geq -cg$
on $[0, \bar \tau]$ for a nonnegative constant $c$. Then \be
\label{firstLestimate}
 %e^{2c\tau}\frac{d^2(p, q, \bar \tau)}{4\tau}+\frac{nc}{3}\tau \geq 
 l(q, \tau) \geq
 e^{-2c\tau}
\frac{d^2(p, q, 0)}{4\tau}-\frac{nc}{3}\tau \ee for each $\tau \in
[0, \bar \tau]$. If we assume instead $Ric \le Cg$ on $[0, \bar
\tau]$ for a nonnegative constant $C$, then \be
\label{secondLestimate} 
%e^{-2C\tau}\frac{d^2(p, q, \bar
%\tau)}{4\tau}-\frac{nC}{3}\tau  \leq 
l(q, \tau) \le
e^{2C\tau}\frac{d^2(p, q, 0)}{4\tau}+\frac{nC}{3}\tau\ee for each
$\tau \in [0, \bar \tau]$.
\end{lem}
\Pf  We first assume a lower bound on the Ricci curvature. By
(\ref{basicboundI}) and (\ref{newL}) we have  for an arbitrary
$\gamma$ from $p$ to $q$ \be \label{unterschrank} {\cal L}(\gamma)
\geq \frac{e^{-2c\tau}}{2}\int_{0}^{\sqrt{\tau}}
|\gamma'|_{g(0)}^2dt-\frac{2nc}{3}\tau^{\frac{3}{2}}\geq
e^{-2c\tau}\frac{d^2(p, q,
0)}{2\sqrt{\tau}}-\frac{2nc}{3}\tau^{\frac{3}{2}}. \ee This leads
to (\ref{firstLestimate}).

%To deduce the first half of  (\ref{firstLestimate}) we apply
%(\ref{basicboundI}) and choose $\gamma$ to be a minimal geodesic
%from $p$ to $q$ with respect to the metric $g(\bar \tau)$. 
The
case of an upper bound for the Ricci curvature is similar, in
which we use (\ref{basicboundII}) instead of (\ref{basicboundI}).
\qed \\

Next we consider Perelman's {\it $\cal{L}$-exponential map}.\\

\noindent {\bf Definition 2}
 The {\it $\cal{L}$-exponential map} $exp^{\cal{L}, \tau}_{p}: T_pM \rightarrow M$ at
time $\tau \in [0, T)$ is defined as follows. For $v \in T_{p}M$,
let $\gamma_v$ denote the $\cal{L}$-geodesic such that
$\gamma_v(0)=p, \lim_{s \rightarrow 0} \sqrt{s} \dot{\gamma}(s)
=v$ (equivalently, $\gamma_v'(0)=2v$). If $\gamma_v$
exists on $[0, \tau]$, we set $exp^{\cal{L},
\tau}_{p}(v)=\gamma_v(\tau).$ Let ${\cal{U}}(\tau)$ denote the
maximal domain of $exp^{\cal{L}, \tau}_{p}$. By
(\ref{newgeodesic}) and basic ODE, ${\cal{U}}(\tau)$ is an open
set and  $exp^{\cal{L}, \tau}_{p}$ is a smooth map from
${\cal{U}}(\tau)$ into $M$.\\

We also have the following extension of the concept of
$\cal{L}$-exponential map. \\

\noindent {\bf Definition 3} For a given reference point $\bar p$
and $0<\varepsilon <\tau$  the ${\cal L}_{\varepsilon,
\tau}$-exponential map $exp^{\cal{L}_{\varepsilon, \tau}}_{\bar
p}$ is defined as follows. For  $v \in T_{\bar p}M$, let
$\gamma_{v,\varepsilon}$ denote the $\cal{L}_{\varepsilon, \tau}$-geodesic such that
$\gamma_{v, \varepsilon}(\varepsilon)=\bar p,  \sqrt{\varepsilon}
\dot{\gamma}(\varepsilon) =v$ (equivalently, $\gamma_{v,
\varepsilon}'=2v$ at $t=\sqrt{\varepsilon}$). If $\gamma_{\varepsilon,v}$ exists
on $[0, \tau]$, we set $exp^{\cal{L}_{\varepsilon, \tau}}_{\bar
p}(v)=\gamma_{\varepsilon, v}(\tau).$

\begin{prop} \label{explemma} Assume that the
sectional  curvature is bounded on $[0, \bar \tau]$ for $\bar \tau \in
(0, T)$.
   Then ${\cal{U}}(\tau)=T_pM$ for each $\tau \in (0, \bar \tau)$.
   A similar statement holds true for $exp^{\cal{L}_{\varepsilon, \tau}}_{\bar
p}$.
\end{prop}
\Pf By the local interior estimates in [S], the sectional curvature
bound on $[0, \bar \tau]$ implies an upper bound on $|\nabla R|$
on $[0, \tau]$ for each $\tau \in (0, \bar \tau)$. Fix $\tau \in
(0, \bar \tau)$ and let $K$ denote an upper bound for 
$|Ric|$ and
$|\nabla R|$ on $M \times [0, \tau]$.

Consider an $\cal L$-geodesic $\gamma$ with initial time $0$,
defined on its maximal interval. We derive from
(\ref{newgeodesic}) \ba {\frac{d}{dt}} |\gamma'|^2&=&
\frac{\partial g}{\partial s}(\gamma', \gamma')\frac{ds} {d t}+2
\gamma' \cdot \nabla_{\frac{d}{dt}} \gamma' = 4t^2 \nabla R \cdot
\gamma'.
%-4tRic(\gamma', \gamma'). 
\ea 
 Consequently, we obtain for
$t \le \sqrt{\tau}$ (as long as $\gamma$ is defined) \be
\label{gammaderivative} |\frac{d}{dt} |\gamma'|^2| \le 
4K t^2 |\gamma'|
\ee
and hence 
\be
|\frac{d}{dt}|\gamma'| | \le 2Kt^2.
\ee
%+ 4Kt |\gamma'|^2 \le 2Kt^3+6Kt|\gamma'|^2 \ee and
%hence 
It follows that 
\be \label{gammabound} |\gamma'| \le
%|\gamma'(0)|^2+\int_0^t 2K u^3e^{-3Ku^2}du) \ee for $t
|\gamma'|(0)+\frac{2}{3}Kt^3.
\ee
%for 
%$ t \le \sqrt{\tau}$. 
%Since $\dot{\gamma}=\gamma'/2\sqrt{s}$ this implies
%\be  \label{gammabound}
%|\dot{\gamma}| \le \frac{1}{\sqrt{s}}e^{2Ks}
%(|v|^2+\int_0^{\sqrt{s}} K u^3e^{-4Ku^2})^{\frac{1}{2}}du.
%\ee
By (\ref{basicboundI})
we then infer 
\ba  \label{gammabound1}
|\gamma'(t)|_{g(0)} \le e^{Kt}(|\gamma'|(0) +\frac{2}{3}Kt^3).
\ea
This gives rise to a uniform  upper bound for the length of $\gamma|_{[0, \tau']}$
for $\tau' \le \tau$ measured in $g(0)$.
By the completeness of $g(0)$ and basic ODE we conclude that
 $\gamma$ is defined on $[0, \tau]$.
\qed  \\

\begin{prop} \label{approachlemma} We have $\cup_{\tau} {\cal{U}}(\tau)=T_pM$. In other
words, the direct limit of ${\cal{U}}(\tau)$ as $\tau \rightarrow
0$ is $T_pM$. Indeed, for each $r>0$, there is $\tau>0$ such that
$B_r(0) \subset {\cal{U}}(\tau)$, where the norm on $T_pM$ is
induced from $g(0)_p$. A similar statement holds true for $
exp^{\cal{L}_{\varepsilon, \tau}}_{\bar p}$.
\end{prop}
\Pf Fix $0<\tau^*<\mbox{min}\{T, 1\}$. Let $r>0$ be given. Let $K$
be an upper bound for $|Ric|$ and $\nabla R$ on $B_{2r}(p, 0)
\times [0, \tau^*]$. Consider $v \in B_r(0)$.
% Applying
%(\ref{gammabound}) and Lemma 2.1 to
By (\ref{gammabound1}) we have for  $\gamma=\gamma_v$
(parameterized in $t$) 
\be |\gamma'(t)|_{g(0)} \le
e^{Kt}(r+\frac{2}{3}Kt^3),  \ee 
as long as $t$ is in the maximal existence interval
of $\gamma$, $t \le \sqrt{\tau^*}$ and $\gamma([0, t]) \subset B_{2r}(p,
0)$.  For such $t$ which also satisfies $t
< \sqrt[3]{\frac{r}{2K}}$ and $t< \frac{1}{K}\ln
\frac{9}{8}$ we then have \be |\gamma'(t')|_{g(0)} <
\frac{3}{2}r \ee for all $t' \in [0, t]$, whence $\int_0^t
|\gamma'|_{g(0)} dt < \frac{3}{2}r$. Consequently, $\gamma([0,
t]) \subset B_{\frac{3}{2}r}(p, 0)$. Since $\gamma(0)=p
$, by continuity we then obtain 
$\gamma([0, t]) \subset B_{\frac{3}{2}r}(p, 0)$ for all $t$ in the
maximal existence interval of $\gamma$ such that $0 \leq t \leq
t^*$, where
$t^*=\mbox{min}\{\sqrt{\tau^*},\sqrt[3]{\frac{r}{2K}},\frac{1}{K}\ln
\frac{9}{8} \}$. This implies in turn that $[0, t^*]$ is contained
in the maximal existence interval of $\gamma$. It follows that $v
\in {\cal{U}}(\tau)$ for each $\tau \in [0, {t^*}^2]$.  \qed

\begin{prop} For each sufficiently small $\tau>0$, $\lexp$ is
a diffeomorphism from a neighborhood of $0$ in $T_pM$ onto a
neighborhood of $p$ in $M$. If $\nabla^2R (p, \tau) \geq 0$ for
each $\tau$, then this holds for each $\tau$. A similar statement
holds for $exp^{\cal{L}_{\varepsilon, \tau}}_{\bar p}$.
\end{prop}
\Pf First note that $0 \in {\cal{U}}(\tau)$. Indeed, the ${\cal
L}$-geodesic $\gamma_0$ is the constant curve $\gamma_0 \equiv p$,
hence it is defined for all $\tau$.  To establish the desired
difeomorphism property, it suffices to show that the differential
of $\lexp$ at $0$ is has zero kernel. For this purpose, consider a
nonzero $v \in T_pM$ and $\lexp(xv)=\gamma_{xv}(\tau)$. Obviously,
${d\lexp}|_0(v)=Y_v(\sqrt{\tau})$, where $Y_v$ (parameterized in
$t$) is the $\cal L$-Jacobi field along $\gamma_0$ associated with
the family of $\lexp$-geodesics $\gamma_{xv}$ (with parameter
$x$). Thus $Y_v(0)=0, \nabla_{\frac{d}{dt}} Y_v(0)=v$. By [(7.7),
P], the $\cal L$-Jacobi equation along an $\cal L$-geodesic
$\gamma$ (parameterized in $s$) is \be \nabla_{\frac{d}{ds}}
\nabla_{\frac{d}{ds}}Y+\frac{1}{2s} \nabla_{\frac{d}{ds}}Y
+Rm(\dot{\gamma}, Y)\dot{\gamma}+2\nabla_YRic(\dot{\gamma}, \cdot)
-\nabla_{\dot{\gamma}}Ric(Y, \cdot)-\frac{1}{2} \nabla^2 R(Y,
\cdot)=0. \ee For $\gamma=\gamma_0$ this becomes, when
parameterized in $t$, \be \frac{d^2Y}{dt^2}-2t^2 \nabla^2R(p,
t^2)(Y, \cdot)=0. \ee It is easy to see that for small $\tau$,
$Y(0)=0$ and $Y(\sqrt{\tau})=0$ imply that $Y \equiv0$. The same
holds for each $\tau$ if $\nabla^2R(p, \tau)\geq 0$ for each
$\tau$. Applying this to $Y_v$ we arrive at the desired
conclusions. \qed

\begin{prop}  \label{expontolemma} If $Ric\geq -cg$ on $[0, \tau]$
for a  nonnegative constant $c$, then there exists a minimal
${\cal{L}}_{0, \tau}$-geodesic from $p$ to $q$ for each $q$.
Consequently, $\lexp$ is onto. 
%The same conclusions hold if $Ric
%\leq Cg$ on $[0, \tau]$ for a nonnegative constant $C$.
\end{prop}
\Pf For a given $q \in M$ we minimize the $\cal{L}$-energy in the
reparametrized form (\ref{newL}) among Sobolev curves which
connect $p$ to  $q$. By the estimate (\ref{firstLestimate}) we can
find a minimizer $\gamma$.  By the standard elliptic regularity,
it is a smooth $\cal{L}$-geodesic connecting $p$ to $q$. Set
$v=\gamma'(0)/2$. Then $\lexp (v)=q$.
\qed  \\

\noindent {\bf Definition 4} 1) We define the   {\it injectivity
domain} $\Omega({\tau})$ at time $\tau$ to be
$$\Omega({\tau})=\{q \in M: \mbox{ there is a unique minimal } {\cal{L}}-\mbox{geodesic }
\gamma: [0, \tau] \rightarrow M
$$
$$\mbox{ with } \gamma(0)=p, \gamma(\tau)=q;
q \mbox{ is not conjugate to } p \mbox{ along }
\gamma\}.$$
Here, ``conjugate" means the same as in ordinary Riemannian geometry of
geodesics, i.e. there is a nontrivial $\cal{L}$-Jacobi field
$J$ along $\gamma$ with $J(0)=0, J(\tau)=0.$

The {\it cut-locus} ${\cal C}({\tau})$ is defined to be
$M-\Omega({\tau})$.

The corresponding concepts, the $\cal{L}_{\varepsilon, \tau}$
injectivity domain $\Omega(\varepsilon, \tau)$ and cut-locus
$\cal{C}(\varepsilon, \tau)$ associated with $\cal{L}_{\varepsilon,
\tau}$-geodesics, are defined in a similar way.

2) The {\it tangential injectivity domain} $\Omega^{T_p}(\tau)$ at time
$\tau$ is defined to be
$$
\Omega^{T_p}(\tau)=\{v \in {\cal{U}}(\tau): \gamma_v|_{[0, \tau]} \mbox{ is a unique minimal }
{\cal{ L}}-\mbox{geodesic s.t. } \tau \mbox{ is not a conjugate}
$$
$$
\mbox{ time}.\}
$$
It is easy to see that $\Omega(\tau)=\lexp(\Omega^{T_p}(\tau))$.

The tangential $\cal{L}_{\varepsilon, \tau}$ injectivity domain
$\Omega^{T_{\bar p}}(\varepsilon, \tau)$ is defined in a similar
way. \\

%\begin{prop} \label{analyticlemma} $\cup_{\tau} C(\tau)$ is a real
%analytic subvariety of dimension $n$ in $M \times [0, T)$. For
%each $\tau$, $C(\tau)$ is a real analytic subvariety of $M$ of
%dimension $n-1$.
%\end{prop}
%\Pf

\begin{lem} \label{smoothlemma} $\cal{C}(\tau)$ is closed in $M$ for each $\tau \in (0, T)$, 
and $\cup_{0<\tau<T} \cal{C}(\tau) \times \{\tau\}$ is closed in $M \times (0,T)$.  Consequently, $\Omega(\tau)$, 
$\Omega^{T_p}(\tau), \cup_{0<\tau<T}\Omega(\tau) \times \{\tau\}$ and $\cup_{0<\tau<T} \Omega^{T_p}(\tau) \times 
\{\tau\}$ 
 are open.  $\lexp$ is a smooth diffeomorphism
from $\Omega^{T_p}(\tau)$ onto $\Omega(\tau)$, depending smoothly
on the parameter $\tau$.   $L(q, \tau)$ is a smooth function on
$\cup_{0<\tau<T}\Omega(\tau) \times \{\tau\}$.

Similar statements hold in the situation of $exp^{\cal{L}_{\varepsilon, \tau}}_{\bar p}$.
In particular, $exp^{\cal{L}_{\varepsilon, \tau}}_{\bar p}$ is a
smooth diffeomorphism from $\Omega^{T_{\bar p}}(, \varepsilon,
\tau)$ onto $\Omega(\varepsilon, \tau)$, depending smoothly on
$\varepsilon$ and $\tau$,  and $L_{\varepsilon, \tau}(q)$ is
smooth on $\cup_{\varepsilon, \tau} \{(\varepsilon, \tau)\} \times
\Omega(\varepsilon, \tau)$, where $L_{\varepsilon, \tau}(q)$ is
defined to be the infimum of $\cal{L}_{\varepsilon, \tau}(\gamma)$
for $\gamma:[\varepsilon, \tau] \rightarrow M$ such that
$\gamma(\varepsilon)=\bar p, \gamma(\tau)=q$.
\end{lem}
\Pf All these can easily be established  by applying the corresponding standard arguments in the theory of ordinary geodesics in
Riemannian geometry.  \qed

\begin{lem} \label{geodesiclemma} Let $\gamma$  be
a minimal ${\cal{L}}_{0, \tau}$-geodesic. Then $\gamma|_{[0,
\tau']} $ is the unique minimal ${\cal{L}}_{0, \tau'}$-geodesic
from $p$ to $\gamma(\tau')$  for any $\tau' \in (0, \tau)$.
Moreover, $\tau'$ is not a conjugate time. Thus $\gamma(\tau') \in
\Omega(\tau')$. We also have $\gamma(\tau') \in \Omega(\varepsilon,
\tau)$ for any $\varepsilon \in (0, \tau)$ and $\tau' \in (\varepsilon, \tau)$, where the reference
point $\bar p$ for $\Omega(\varepsilon, \tau)$ is chosen to be
$\gamma(\varepsilon)$.

 As a consequence, we have
$\Omega^{T_p}(\tau_2) \subset \Omega^{T_p}(\tau_1)$ for $\tau_2
\geq \tau_1$.
\end{lem}
\Pf The arguments in the theory of ordinary geodesics can be applied
directly.
\qed  \\

By the first variation formula [(7.1), P1], if $q \in \Omega(\tau)$ and $\gamma$ is the unique minimal
$\cal{L}$-geodesic from $p$ to $q$, then  we have
\be \label{gradientformula}
\dot{\gamma}(s)=\nabla l(\gamma(s)), \gamma'(s(t))=\nabla L(\gamma(s(t)))
\ee
for $s \in [0, \tau]$ and $s(t)=t^2$.

\begin{prop} \label{lipschitzlemmaI} Let $\bar \tau \in (0, T)$.  Assume that $Ric\geq -cg$
on $[0, \bar \tau]$ for a nonnegative constant $c$. Then $L(\cdot,
\tau)$ is locally Lipschitz with respect to the metric $ g(\tau)$
for each $\tau \in (0, \bar \tau]$. Moreover, for each compact
subset $E$ of $M$, there are positive constants $A_1$ and $A_2$
such that $\sqrt{\tau}L\le A_1$  on $E \times (0, \bar \tau]$ and
\be \label{gammaderivativebound} |\dot{\gamma}(s)|^2\le
\frac{A_2}{s}(1+\frac{1}{\tau}) \ee for $s \in (0, \tau]$, where
$\tau \in (0, \bar \tau]$ and $\gamma$ denotes an arbitrary
minimal ${\cal{L}}_{0, \tau}$-geodesic from $p$ to $q$ for $q \in
E$.
\end{prop}
\Pf We first derive an upper bound for $\sqrt{\tau} L(q, \tau)$ on
$B_{\rho}(p, \bar \tau) \times (0, \bar \tau]$ for a given
$\rho>0$.  By smoothness, there is a positive constant $C$ such
that $R \leq C$ on $B_{\rho}(p, \bar \tau) \times [0, \bar \tau]$.
For $q \in B_{\rho}(p, \bar \tau)$ and $\tau \in (0, \bar \tau]$
we choose a minimal geodesic $\gamma: [0, \sqrt{\tau}] \rightarrow
B_{\rho}(p, \bar \tau)$ from $p$ to $q$ with respect to $g(\bar
\tau)$. By (\ref{basicboundI}) and (\ref{newL})  we have
$$
{\cal{L}}(\gamma) \leq \int_0^{\sqrt{\tau}} (e^{2c(\bar \tau-t)}
|\gamma'|^2_{g(\bar \tau)}+2Ct^2)dt \leq e^{2c\bar \tau}
 \frac{d(p, q, \bar \tau)}{\sqrt{\tau}}+3C \tau^{\frac{3}{2}}.
 $$
 It follows that
 $$
 \sqrt{\tau}L(q, \tau) \leq A({\rho}),
 $$
 where $A(\rho)= e^{2c\bar \tau} \rho + 3C \bar \tau^2$.

Next consider a given $\rho>0$. Choose $\rho_1$ such that
$B_{\rho}(p, 0) \subset B_{\rho_1}(p, \bar \tau)$.
%By smoothness, the Ricci curvature is bounded from above on
%$B_{\rho}(p, 0) \times [0, \bar \tau]$. By the arguments in the
%proof of Lemma \ref{firstllemma} we then obtain an upper bound
%$\hat L$ for $\sqrt{\tau}|L(q, \tau)|$ on $B_{\rho}(p, 0) \times
%(0, \bar \tau]$.
We set $\rho^{*}=\max\{e^{c\bar \tau} \sqrt{\frac{4nc}{3}{\bar
\tau}^2 +2A(\rho_1)}, 2\rho\}$.
% Fix a $\tau^* \in (\bar \tau,T)$.
By the smoothness of $g$, there is an upper bound $K$ for $|Ric|$
and $|\nabla R|$
%the sectional curvature on
%$B_{2\rho^{*}}(p, 0) \times [0, \tau^{*}]$ is bounded. By the
%local interior estimates in [S], we then obtain an upper bound $K$
%for $|\nabla R|$
on $B_{\rho^{*}}(p, 0) \times [0, \bar \tau]$.
%We choose $K$ such that it is also an upper bound for $|Ric|$ on this region.

Now consider $q_1, q_2 \in B_{\rho}(p, 0)$ and $\tau \in (0, \bar
\tau]$. Let $\gamma_i$ be a minimal ${\cal{L}}_{0, \tau}$-geodesic
from $p$ to $q_i, i=1, 2.$ (By Lemma \ref{expontolemma}, they
exist.) Let $\gamma_0: [0, 1] \rightarrow M$ be a minimal geodesic
from $q_1$ to $q_2$ with respect to $g(0)$.  By the choice of
$\rho^*$, the image of $\gamma_0$ is obviously contained in
$B_{\rho^{*}}(p, 0)$. We claim that the images of $\gamma_1$ and
$\gamma_2$ are also contained in $B_{\rho^*}(p, 0)$.  Indeed, we
have
\be \int_0^{\tau'} \sqrt{s}(R+|\dot{\gamma_i}|^2)ds \le
L(q_i, \tau) -\int_{\tau'}^{\tau} \sqrt{s} Rds \le
\frac{A(\rho_1)}{\sqrt{\tau}}+ \frac{2nc}{3} \tau^{\frac{3}{2}}
\ee
 for $i=1, 2$ and $\tau'
\in [0, \tau]$.  By (\ref{geodesiclemma}) and
(\ref{firstLestimate}) (applied to $\tau'$) we then deduce \be
d(p, \gamma_i(\tau'),0)^2 \leq
4e^{2c\tau'}(\frac{1}{2}\sqrt{\tau'} L(q_i, \tau') + \frac{nc}{3}
\tau'^2) \leq \rho^{*2} \ee for $i=1, 2$ and $\tau' \in [0,
\tau]$. It follows that
 the images of
$\gamma_1$ and $\gamma_2$ are contained in $B_{\rho^{*}}(p, 0)$.

Next we estimate $|\dot{\gamma_1}|$ and $|\dot{\gamma_2}|$. It is
more convenient to handle $\gamma_1'$ and $\gamma_2'$. By the
arguments in the proof of Lemma \ref{explemma} we deduce 
\be
|\gamma_i'|(t_2) 
\ge |\gamma_i'|(t_1)-\frac{2}{3}K\bar \tau^3 \ee for
$i=1, 2$ and $t_1, t_2 \in [0, \sqrt{\tau}]$ and hence \be
\label{tangentbound} 
% \frac{1}{2}|{\gamma_i}'(t_1)|^2-C_1 \le
|{\gamma_i}'(t_2)|^2\ge  \frac{1}{2}| {\gamma_i}'(t_1)|^2-C\ee for
$t_1, t_2 \in [0, \sqrt{\tau}], i=1, 2$, where 
$C=\frac{4}{9}K^2\bar \tau^6$.  It follows
that \be \label{gammaprime} 4s |\dot{\gamma_i(s)}|^2
=|\gamma_i'|^2 \le 8(\frac{L(q_i,
\tau)}{\sqrt{\tau}}+\frac{2nc}{3}\tau) +2C \ee for $i=1,
2$ and $s \in [0, \tau]$.

To proceed, we set $d=d(q_1, q_2, 0)$ and assume that $d
<\frac{\tau}{4}$. We define $\hat \gamma_1(s)=\gamma_1(s)$ for $ s
\in [0, \tau-2d] $, $\hat
\gamma_1(s)=\gamma_1(\tau-2d+2(s-\tau+2d))$ for $s\in [\tau-2d,
\tau-d]$ and $\hat \gamma_1(s)=\gamma_0(\frac{1}{d}(s-\tau+d))$
for $s \in [\tau-d, \tau]$. Then we have \ba L(q_2, \tau)&\le&
{\cal{L}}(\hat \gamma_1)\le L(q_1, \tau)
-\int_{\tau-2d}^{\tau}\sqrt{s}R(\gamma_1)ds+\int^{\tau-d}_{\tau-2d}\sqrt{s}
(R(\gamma_1)+4|\dot{\gamma_1}|^2)ds \nonumber \\
&&+\int_{\tau-d}^{\tau}\sqrt{s}(R( \gamma_0)+
\frac{1}{d^2}|\dot{\gamma_0}|^2)ds, \ea where the arguments for
$\gamma_1$ and $\gamma_0$ in the second and third integrals
correspond to the defintion of $\hat \gamma_1$.  We have \be
-\int_{\tau-2d}^{\tau} \sqrt{s} R(\gamma_1)ds \leq \frac{2nc}{3}
(\tau^{\frac{3}{2}}-(\tau-2d)^{\frac{3}{2}}),  \ee \be
\int_{\tau-2d}^{\tau-d} \sqrt{s}R(\gamma_1) \leq \frac{2}{3}nK
((\tau-d)^{\frac{3}{2}}-(\tau-2d)^{\frac{3}{2}}), \ee
and
 \be \int_{\tau-d}^{\tau}\sqrt{s}R(
\gamma_0)ds \leq
\frac{2}{3}nK(\tau^{\frac{3}{2}}-(\tau-d)^{\frac{3}{2}}). \ee By
(\ref{basicboundII}) we have $|\dot{\gamma_0}|^2\le
e^{2K\tau}d^2$, hence \be
\int_{\tau-d}^{\tau}\sqrt{s}\frac{1}{d^2}|\dot{\gamma_0}|^2ds \leq
\frac{2}{3}e^{2K\tau}(\tau^{\frac{3}{2}}-(\tau-d)^{\frac{3}{2}}).
\ee On the other hand, we infer from (\ref{gammaprime}) that \be
\int_{\tau-2d}^{\tau-d}4\sqrt{s}|\dot{\gamma_1}|^2 ds \leq
4(4\frac{A(\rho_1)}{{\tau}}+\frac{8nc}{3}\tau
+C)((\tau-d)^{\frac{1}{2}}-(\tau-2d)^{\frac{1}{2}}). \ee We deduce \be \label{oben} L(q_2, \tau)\le
L(q_1, \tau)+I(\tau, d), \ee where   $$ I(\tau, d)  =
\frac{2}{3}(2nc+nK+
e^{2K\tau})(\tau^{\frac{3}{2}}-(\tau-2d)^{\frac{3}{2}})
$$
\ba +4(4\frac{A(\rho_1)}{{\tau}}+\frac{8nc}{3}\tau
+C)((\tau-d)^{\frac{1}{2}}-(\tau-2d)^{\frac{1}{2}}). \ea
 Similarly, we have \be \label{unten} L(q_1, \tau) \le L(q_2,
\tau)+I(\tau, d). \ee The desired Lipschitz continuity follows.
The estimate (\ref{gammaderivativebound}) follows from
(\ref{gammaprime}).

Finally, we would like to point out that the local Lipschitz
continuity of $L(\cdot, \tau)$ also follows from its local
semiconcavity, which is given by Lemma \ref{concavelemma} below.
Note however that the proof of Lemma \ref{concavelemma} below uses
some arguments here. \qed

\begin{prop} \label{lipschitzlemmaII}
Assume that the Ricci curvature is bounded from below on
$[0, \bar \tau]$. Then $L(q, \cdot)$ is locally Lipschitz
on $(0, \bar \tau]$ for every $q\in M$. Moreover, ${\tau}^{\frac{3}{2}}|L_{\tau}|$
is bounded on $E \times (0, \bar \tau]$ for each compact subset $E$ of
$M$.
\end{prop}
\Pf This is similar to the proof of Proposition
\ref{lipschitzlemmaI} above. Fix $\rho>0$ and let $ \rho^{*}$ and
$K$ have the same meanings as in the proof of Proposition
\ref{lipschitzlemmaI}.
 Consider $q \in B_{\rho}(p, 0)$ and
$\tau_1, \tau_2 \in (0, \tau]$ such that $ \tau_1<\tau_2$ and
$\tau_2 < 2\tau_1$. Choose a minimal ${\cal{L}}_{0,
\tau_1}$-geodesic $\gamma_1$ from $p$ to $q$ and a minimal
${\cal{L}}_{0,  \tau_2}$-geodesic $\gamma_2$ from $p$ to $q$. As
in the proof of Proposition \ref{lipschitzlemmaI}, the images of
$\gamma_1$ and $\gamma_2$ are contained in $B_{\rho^{*}}(p, 0)$.
We define $\hat \gamma_1(s)=\gamma_{1}(s), s \in [0, \tau_1]$ and
$\hat \gamma_{1}(s)=q, s \in [\tau_1, \tau_2]$. Then \be
\label{tauestimateI} L(q, \tau_2)\leq {\cal{L}}_{0, \tau_2}(\hat
\gamma_{1})\le L(q, \tau_1)
+\int_{\tau_1}^{\tau_2}\sqrt{s}R(q, s)ds\\
\leq L(q, \tau_1)+\frac{2}{3}nK(\tau_2^{3/2}-\tau_1^{3/2}).
\ee

Next we set $\tau_3 =2\tau_1-\tau_2$, $\hat \gamma_{2}(s)=\gamma_{2}(s)$ for $s \in
[0, \tau_3]$
and
$\hat \gamma_{2}(s)=\gamma_{2}(\tau_3+2(s-\tau_3))$ for
$s \in [\tau_3, \tau_1]$. Then
\ba \label{tauestimateII}
L(q, \tau_1)&\le& {\cal{L}}_{0, \tau_1}(\hat \gamma_2) \le
L(q, \tau_2)-\int_{\tau_3}^{\tau_2} \sqrt{s}R(\gamma_2)ds \nonumber\\ &&
+\int_{\tau_3}^{\tau_1}\sqrt{s}(R(\gamma_{2})+
4|\dot{\gamma_2}|^2)ds,
\ea
where the argument of $\gamma_2$ in the last integral on the right hand side
is $\tau_3
+2(s- \tau_3)$. Applying (\ref{gammaderivativebound}) we then obtain
\be \label{tauestimateIII}
L(q, \tau_1)\le L(q, \tau_2)+ \frac{2n(c+K)}{3}(\tau_1^{\frac{3}{2}}-\tau_3^{\frac{3}{2}})
+8A_2(1+\frac{1}{\tau})(\tau_1^{\frac{1}{2}}-\tau_3^{\frac{1}{2}}).
\ee
Clearly, (\ref{tauestimateI}) and (\ref{tauestimateIII}) imply the desired
Lipschitz continuity and derivative bound.
\qed

\begin{prop} \label{lipschitzlemmaIII} Assume that the Ricci curvature is bounded from below on
$[0, \bar \tau]$. Then $L$ is a locally Lipschitz function on $M \times (0, \bar \tau]$. 
\end{prop}
\Pf  Combine Proposition \ref{lipschitzlemmaI} and Proposition \ref{lipschitzlemmaII}.
To be more precise, we have $|L(q_1,\tau_1)-L(q_2, \tau_2)| \le |L(q_1, \tau_1)-L(q_1, \tau_2)|
+|L(q_1, \tau_2)-L(q_2,\tau_2)|$. We apply the above two propositions to handle the 
two terms on the right hand side to obtain the desired Lipschitz bound.  \qed

\begin{prop} \label{concavelemma} Assume that the Ricci curvature
is bounded from below on $[0, \bar \tau]$. Then $l(\cdot, \tau)$
is locally semi-concave for each $\tau \in (0, \bar \tau]$, i.e.
for every point $q \in M$ there is a smooth function $\phi$ on a
neighborhood $U_q$ of $q$ such that $l(\cdot, \tau)+\phi$ is
concave in the sense that the composition of $l(\cdot, \tau)+\phi$
with every geodesic in $U_q$ is a concave function.
\end{prop}
\Pf   By [(7.9), P] we have for each $\tau \in (0, T), q \in
\Omega({\tau}) $  and $v \in T_q M$ \be  \label{hessian} Hess_L(v,
v)\leq \frac{1}{\sqrt{\tau}}|v|^2-2\sqrt{\tau}Ric(v, v)-
\int_0^{\tau}\sqrt{s}H(X, Y)ds, \ee where $X=\dot{\gamma}$ with
$\gamma$ denoting the unique minimal $\cal{L}$-geodesic from $p$
to $q$, $Y$ is a suitable extension of $v$ along $\gamma$ such
that $|Y(s)|^2=\frac{s}{\tau}|v|^2$, and
$$
H(X, Y)=-\nabla_Y\nabla_YR-2<Rm(Y, X)Y, X>-4(\nabla_XRic(Y,
Y)-\nabla_YRic(Y, X))
$$
\be  \label{harnack} -2Ric_{\tau}(Y, Y)+2|Ric(Y,
\cdot)|^2-\frac{1}{s}Ric(Y, Y). \ee

To estimate $H(X, Y)$ we fix $\rho>0$ and assume $q \in
B_{\rho}(p, 0)\cap \Omega(\tau)$ and $\tau \in (0, \bar \tau]$.
Let $\rho^{*}$ be given in the proof of Proposition
\ref{lipschitzlemmaI}. As in the proof of Proposition
\ref{lipschitzlemmaI}, the smoothness of $g$ implies an upper
bound $C$ for $|\nabla^2 R|, |\nabla Ric|, |Ric_{\tau}|, |Rm|$ and
$|Ric|$ on $B_{\rho^*}(p, 0) \times [0, \bar \tau]$.  By the proof
of Proposition \ref{lipschitzlemmaI}, $\gamma$ is contained in
$B_{\rho^*}(p, 0)$. Hence we have $|H(X, Y)| \le
\frac{s}{\tau}(C(3+2C+\frac{1}{s})+8C|X|+2C|X|^2)|v|^2.$ Applying
(\ref{gammaderivativebound}) we then deduce  $|H(X, Y)| \le
\frac{C_1}{\tau}(s+1+\frac{1}{\tau})|v|^2$ for a positive constant
$C_1$. It follows that \be  \label{hestimate} Hess_L(v, v)\le C_2
|v|^2 \ee for a positive constant $C_2=C_2(\bar \tau)$. 
(Note that
if  the curvature operator is nonnegative, then $H(X, Y)$ can be
estimated as in [7.2, P].) 
% namely Hamilton's Harnack inequality
%implies $H(X, Y)\geq -R(\frac{1}{s}+ \frac{1}{\tau_0-s})|Y|^2$,
%where $\tau <\tau_0<T$.
 %We then deduce from (\ref{hessian})
%\be \label{hestimate}
%Hess_L(Y, Y) \leq \frac{C}{\sqrt{\tau}}|Y|^2,
%\ee
%where $C$ depends on $\tau_0$
%One can e.g. choose $\tau_0=\frac{1}{2}(\tau +T)$ if $T$ is
%finite.)

We claim that (\ref{hestimate}) holds true for all $q \in
B_{\rho}(p, 0)$ in the sense of barriers, provided that $C_2$ is
chosen large enough. This means that for each point $q \in
B_{\rho}(p, 0)$ and each $\varepsilon>0$ we can find a smooth
function $f$ on a neighborhood of $q$ (called an {\it
$\varepsilon$-barrier} at $q$) such that $f\geq L(\cdot, \tau),
f(q)=L(q, \tau)$ and $Hess_{f}(q)(v, v)\leq (C_2
+\varepsilon)|v|^2$. Consider $q \in B_{\rho}(p, 0)$. (We can
assume that $q \in {\cal{C}}(\tau)$.) Choose a minimal $\cal{L}$-geodesic
$\gamma$ from $p$ to $q$. For a given $\varepsilon>0$ we define
\be f=L(\gamma(\varepsilon), \tau)+L_{\varepsilon, \tau}(q), \ee
where $L_{\varepsilon, \tau}$ is defined in Lemma
\ref{smoothlemma} with the reference point $\bar
p=\gamma(\varepsilon)$.
% Note that $d(\gamma(\varepsilon), q) \le
%\rho+\rho^*$.
By Lemma \ref{smoothlemma}, $L_{\varepsilon, \tau}$
is smooth at $q$. We can estimate its Hessian at $q$ in the same
fashion as above. Indeed, all the relevant lemmas can easily be
extended to the situation of $L_{\varepsilon, \tau}$. Then one
infers readily that $f$  is an $\varepsilon$-barrier at $q$.

For each $q \in M$ we choose a suitable smooth function on a
neighborhood of $q$ (for example $\phi=-C'd(q, \cdot, \tau)^2$ for
a suitable $C'$) and deduce that \be Hess_{L+\phi}\leq 0 \ee on a
neighborhood of $q$ in the sense of barriers. The maximum
principle then implies that $L+\phi$ is concave in this
neighborhood (see e.g. [Y3]). \qed

\begin{lem} \label{cutlocuslemma} Assume that the Ricci curvature is
bounded from below on $[0, \tau]$ for $\tau \in (0, T)$. Then the
cut-locus $\cal{C}(\tau)$ is a closed set of measure zero in $M$. Consequently,
$\cup_{0<\tau<T} \cal{C}(\tau) \times \{\tau\}$ is a closed set of measure zero
in $M \times (0, T)$, 
provided that the Ricci curvature is bounded from below on $[0,
\tau]$ for each $\tau \in [0, T)$. 
\end{lem}
\Pf  Set ${\cal{B}}(\tau)=\{q \in M: \exists \mbox{ more than one minimal
} {\cal{L}}_{0, \tau}-\mbox{geodesics from } p \mbox{ to } q \}$
and ${\cal{D}}(\tau)=\{ q \in M: \exists \mbox{ a unique minimal }
{\cal{L}}_{0, \tau}-\mbox{geodesic } \gamma \mbox{ from } p \mbox{
to } q, q \mbox{ is conjugate }$ to $ p$ along  $\gamma\}. $ By
Lemma \ref{expontolemma}, we have $\cal{C}({\tau})=\cal{B}({\tau}) \cup
\cal{D}({\tau})$.  As in the theory of ordinary geodesics, ${\cal{D}}(\tau)$ is
contained in the set of critical values of $\lexp$. By Sards'
theorem, it has zero measure. On the other hand, $L(\cdot, \tau)$
is obviously non-differentiable at any point of ${\cal{B}}({\tau})$. (We 
would like to thank G.~Wei for helpful discussions on this 
point.) Since
$L(\cdot, \tau)$ is almost everywhere differentiable by
Proposition \ref{lipschitzlemmaI}, $\cal{B}({\tau})$ has zero measure.
 It follows that $\cal{C}({\tau})$ has zero measure.

Next we assume that the Ricci curvature is bounded below on $[0,
\tau]$ for each $\tau$. By Lemma \ref{smoothlemma}, $\cup_{0<\tau<T}\cal{C}(\tau) \times 
\{\tau\}
$ is closed in $M \times (0,T)$ and hence measurable. Then the Fubini theorem implies that 
it has measure zero.

Instead of using the Lipschitz property of $l$ in the above proof, we can also use an idea 
suggested in [KL]. 
% to deduce $\cal{C}({\tau})=\cal{B}({\tau}) \cup
%\cal{D}({\tau})$. 
By Sards' theorem,
we only need to show that $\cal{B}(\tau)^*$ has zero measure, where
$\cal{B}(\tau)^*$ is the intersection of $\cal{B}(\tau)$ with the set of
regular   values of $\lexp$. Consider $q \in \cal{B}(\tau)^*$. Then
there are $v_1, v_2 \in T_pM$ such that $v_1 \not = v_2,
\lexp(v_1)=\lexp(v_2)=q,$ and $L(v_1, \tau)=L(v_2, \tau)$, where
$L(v, \tau)={\cal{L}}(\gamma_v)$. Since $v_1$ and $v_2$ are
non-critical for $\lexp$, there are disjoint neighborhoods $U_1$
of $v_1$ and $U_2$ of $v_2$ such that $F_1=\lexp|_{U_1}$ and
$F_2=\lexp|_{U_2}$ are diffeomorphisms onto their common image
$U$, which is a neighborhood of $q$.

To proceed, we define $L_*(v, w)=L(v, \tau)-L(w, \tau)$, and set
$S=\{(v, w) \in U_1 \times U_2: F_1(v)=F_2(w)\}$. Obviously, $S$
is an $n$-dimensional submanifold of $U_1 \times U_2$. We claim
that $0$ is a regular value of  $L_*|_S$. Indeed, consider a curve
$(v(t), w(t))$ in $S$ which represents a tangent vector $(v'(0),
w'(0))$ of $S$ at a given point $(v(0), w(0))$. Since
$\lexp(v(t))=\lexp(w(t))$, we have \be \label{dequal}
d(\lexp)_{v(0)}(v'(0))=d(\lexp)_{w(0)}(w'(0)). \ee On the other
hand, by the first variation formula [(7.1), P] for the $\cal{L}$
energy, we have \be \label{lequal} \frac{d L^*(v(t),
w(t))}{dt}(0)= 2\sqrt{\tau}(< \dot{\gamma}_{v(0)}(0),
Y_1>-<\dot{\gamma}_{w(0)}(0), Y_2>), \ee where $Y_1=
d(\lexp)_{v(0)}(v'(0))$ and $Y_2=d(\lexp)_{w(0)}(w'(0))$. Since
$v(0) \not = w(0)$, we have $\dot{\gamma}_{v(0)}(\tau) \not =
\dot{\gamma}_{w(0)}(\tau)$. It follows that $dL^*_{((v(0), w(0))}
((v'(0), w'(0))) \not =0$. By the implicit function theorem,
$L_*|_S^{-1}(0)$ is an $(n-1)$-dimensional submanifold of $U_1
\times U_2$. Consequently, $S^*=\pi_1 \circ F(L^*|_S^{-1}(0))$ is
an $(n-1)$-dimensional submanifold, where $F=(F_1, F_2)$ and
$\pi_1$ denotes the projection from $U \times U$ to the first
factor. We call $S^*$ a {\it local container} for $\cal{B}(\tau)^*$.

It is easy to see that $\cal{B}(\tau)^*$ is contained in a countable
union of local containers. Hence it has zero measure.

An alternative argument  was suggested by Perelman. By Proposition
\ref{concavelemma} and Aleksandrov's theorem (see [Y3]), $L(\cdot,
\tau)$ is twice differentiable almost everywhere. Consequently,
 $\cal{B}(\tau)$ has measure zero. (On the other hand, one
can show that at a point in $\cal{D}(\tau)$, $L(\cdot, \tau)$ cannot be
twice differentiable. This also implies that ${\cal{D}}(\tau)$ has measure zero and 
hence can substitute for the use of Sards' theorem.) 
%(Of course, Aleksandrov's theorem yields more information than the
%mere statement that ${\cal{C}}(\tau)$ has measure zero.)

\qed 

\begin{lem} \label{measurablelemma} Assume that the Ricci curvature is 
bounded from below on $[0, \tau]$ for each $0<\tau<T$. Then $\nabla l$ and 
$l_{\tau}$ exist almost everywhere and are measurable on $M \times (0,T)$.
\end{lem} 
\Pf This follows from Lemma \ref{smoothlemma} and Lemma \ref{cutlocuslemma}, 
or from Proposition \ref{lipschitzlemmaIII}.

\begin{theo} \label{3estimateslemma}
Assume that the curvature operator is
nonnegative for each $\tau \in [0, T)$.  For each $\bar \tau \in (0, T)$ there is a positive 
constant depending only on the dimension $n$ and the magnitude of $\frac{\bar \tau}{T-\bar \tau}$ such that 
\be \label{Restimate}
R\le \frac{Cl}{\tau}
\ee
everywhere on $M \times (0, \bar \tau]$, 
\be \label{nablaestimate}
|\nabla l|^2 \le \frac{Cl}{\tau}
\ee
almost everywhere in $M$ for each $\tau \in (0, \bar \tau]$,
\be \label{lipschitzestimate}
|\sqrt{l}(q_1, \tau)-\sqrt{l}(q_2, \tau)|\le \sqrt{\frac{C}{4\tau}}d(q_1, q_2, \tau)
\ee
for all $\tau \in (0, \bar \tau]$ and all $q_1, q_2 \in M$, 
and
\be  \label{tauestimate}
|l_{\tau}| \le \frac{Cl}{\tau}
\ee
almost everywhere in $(0, \bar \tau]$ for each $q \in M$.
(Note that $\frac{\bar \tau}{T-\bar \tau}$ is understood to be 
zero when $T=\infty$. Thus $C$ depends only on $n$ in this case.)
Moreover, we have the following  Harnack inequality 
\be \label{tauharnack}
(\frac{\tau_1}{\tau_2})^C \le \frac{l(q, \tau_2)}{l(q, \tau_1)}\leq
(\frac{\tau_2}{\tau_1})^C
\ee
for all $q \in M$ and $\tau_1, \tau_2 \in (0, \bar \tau]$ with $\tau_2>\tau_1$.
\end{theo}
\Pf By [(7.16), P] we have for each $\bar \tau \in (0, T)$
\be  \label{nabla-l} |\nabla l|^2 +R
\le C {\frac{l}{\tau}} \ee on $\cup_{0<\tau \le \bar \tau}\Omega(\tau) \times
\{\tau\}$ for a positive constant $C$ depending only on the dimension $n$ and the magnitude of 
$\frac{\bar \tau}{T-\bar \tau}$. The estimates (\ref{Restimate}) and (\ref{nablaestimate}) follow from this,
Proposition \ref{lipschitzlemmaI}, 
and Lemma \ref{cutlocuslemma}.  
Now the estimate (\ref{nablaestimate}) can be rewritten as 
\ba \label{rootestimate}
|\nabla \sqrt{l}|^2 \le \frac{C}{4\tau},
\ea
which implies  (and is equivalent to) (\ref{lipschitzestimate}). Indeed, given $\tau \in (0, \bar \tau]$ and $q \in 
M$, we can apply (\ref{rootestimate}) to derive (\ref{lipschitzestimate}) along almost every radial geodesic 
(for the metric $g(\tau)$)
starting at $q$. By continuity, it holds along every radial geodesic. Hence (\ref{lipschitzestimate}) holds 
for all $q_1, q_2 \in M$.
  
Next we derive (\ref{tauestimate}). Fix $q \in M$. By Propsition \ref{lipschitzlemmaII}, 
$l_{\tau}(q, \tau)$ exists for almost everywhere $\tau$.  Consider $\tau \in (0, \bar \tau]$ such that 
$l_{\tau}(q, \tau)$ exists.  Observe that 
(\ref{tauestimate}) is invariant
under the rescaling $g(\tau) \rightarrow \frac{1}{a} g(a \tau)$,
hence it suffices to prove it in the case that $\bar \tau>1$ and $\tau=1$.  
We consider $\tau_1=1, \tau_2 \in (1, \bar \tau]$ with $\tau_2<2$
and the curves  $\gamma_1, \gamma_2$ and $\hat \gamma_2$
as in the proof of Proposition \ref{lipschitzlemmaII}.
 By Hamilton's Harnack
inequality ([(11.1), P]), $R_{\tau} \le 0$ and hence $R(q, s) \le
R(q, 1)$  for $s \ge 1$. By
(\ref{Restimate}) we then infer as in (\ref{tauestimateI}) \be \label{newtauI} L(q, \tau_2)
\le L(q, 1)+C L(q, 1)({\tau_2}^{\frac{3}{2}}-1). \ee On the other
hand, by Lemma \ref{geodesiclemma} we have $\gamma_2(s) \in \Omega(s)$ for 
$s \in (0, \tau_2)$.  Hence we can apply (\ref{nabla-l}) in
(\ref{tauestimateII}) to deduce \be \label{newtauII} L(q, 1) \le
L(q, \tau_2)+CL(q, \tau_2)(\sqrt{\tau_2}-\sqrt{2-\tau_2}). \ee
Obviously, (\ref{newtauI}) and (\ref{newtauII}) imply
$|l_{\tau}(q, 1)| \le Cl(q, 1)$  for a positive constant $C$ depending only on 
$n$.

Integrating (\ref{tauestimate}) yields the Harnack estimate (\ref{tauharnack}).
\qed \\

%Alternatively, we have by [(7.5), P] and [(7.6), P]
%\be \label{tauderivative}
%L_{\tau}=\sqrt{\tau}R-\frac{|\nabla L|^2}{4\sqrt{\tau}}
%\ee
%on $\cup_{\tau} \Omega(\tau)$.
%Applying Lemma \ref{cutlocuslemma} and (\ref{nablaestimate}) we also obtain
%(\ref{tauestimate}). \qed \\

Similar estimates for $l$ hold in the case of bounded  sectional curvature. 

\begin{prop} \label{boundedcurvaturelemma} Assume that the sectional curvature
is bounded on $[0, \bar \tau]$. Then there is a positive constant
$C=C(\tau^{*})$ for every $\tau^{*}
\in (0, \bar \tau)$ with the following properties. For each $\tau
\in (0, \tau^{*}]$ we have
\be \label{newnablaestimate}
|\nabla l|^2 \le \frac{C}{\tau}(l+\tau+1)
\ee
almost everywhere in $M$. For each $q \in M$ we have
\be \label{newtauestimate}
|l_{\tau}| \le \frac{C}{\tau}(l+\tau+1)
\ee
almost everywhere in $(0, \tau^{*}]$.
\end{prop}
\Pf Consider $\tau^{*} \in (0, \bar \tau), \tau \in (0, \tau^{*}]$
and $q \in M$. By the assumption and the local interior estimates in [S], we have global bounds for $|Rm|$ and $|\nabla R|$ 
on $[0, \tau]^*$.  By the arguments in the proof of Proposition
\ref{lipschitzlemmaI} we then deduce for a minimal ${\cal{L}}_{0,
\tau}$-geodesic  $\gamma$ from $p$ to $q$ \be
\label{newgammaderivativebound} |\dot{\gamma}|^2 \le
\frac{C}{s}(l(q, \tau)+\tau+1) \ee for a positive constant
$C=C(\tau^{*})$.  Taking $s=\tau$ in
(\ref{newgammaderivativebound}) and applying Lemma
\ref{cutlocuslemma} and (\ref{gradientformula}) we then arrive at (\ref{newnablaestimate}).

The estimate (\ref{newtauestimate}) follows from
(\ref{newgammaderivativebound}) and the arguments in the proof of
Proposition \ref{lipschitzlemmaII}. \qed

\begin{lem} \label{laplacelemma} Assume that
the Ricci curvature is bounded from below on $[0, \bar \tau]$.
Then there holds for every $\tau \in (0, \bar \tau]$
 \be
\label{laplaceI} \int_M l \Delta \phi dq \leq \int_{*M} \phi
\Delta l dq
\ee
for nonnegative smooth functions $\phi$ with
compact support, where the integral $\int_{*M}$ means
$\liminf_{\epsilon \rightarrow 0} \int_{M- U_{\epsilon}}$, with
$U_{\epsilon}=U_{\epsilon}({\cal{C}}(\tau))$ denoting the
$\epsilon$-neighborhood of ${\cal{C}}(\tau)$ ($\epsilon >0$).
Consequently, we have \be  \label{weaklaplace} -\int_M \nabla l
\cdot \nabla \phi dq \leq \int_{*M} \phi \Delta l dq \ee for
nonnegative Lipschitz functions $\phi$ with compact support.
\end{lem}
\Pf Consider $q_0 \in M$. By Proposition~\ref{concavelemma}, there
is a neighborhood $U$ of $q_0$ and a smooth function $\psi$ on $U$
such that $l+\psi$ is concave. We can assume that $l+\psi$ is
actually strictly concave, i.e. it is the sum of a concave
function and a smooth concave function with negative Hessian. By
[GW] or [Y3], there exists a sequence of smooth concave functions
$f_k$ with negative Hessian on a neighborhood $\hat U \subset U$
of $q_0$ such that: 1) $f_k$ converge uniformly to $l+\psi$ on
$\hat U$, and 2) the derivatives of $f_k$ converge uniformly to
the derivatives of $l+\psi$ on $\hat U -U_{\epsilon}({\cal{C}}({\tau}))$
for each $\epsilon>0$.

Let $\phi$ be a nonnegative smooth function with compact support
contained in $\hat U$. Setting
$U_{\epsilon}=U_{\epsilon}({\cal{C}}({\tau}))$ we then have \be \int_M f_k
\Delta \phi dq =\int_M \Delta f_k \phi dq= \int_{M-U_{\epsilon}}
\Delta f_k \phi dq+\int_{U_{\epsilon}} \Delta f_k \phi dq\leq
\int_{M-U_{\epsilon}} \Delta f_k \phi dq. \ee Taking limit we
deduce \be \int_M(l+\psi) \Delta \phi dq \leq
\int_{M-U_{\epsilon}} \Delta (l+\psi) \phi dq. \ee It follows that
\be \int_M (l+\psi)\Delta \phi dq \leq \liminf_{\epsilon
\rightarrow 0} \int_{M-U_{\epsilon}} \Delta (l+\psi) \phi dq. \ee
Since ${\cal{C}}({\tau})$ is closed and has zero measure by Lemma
\ref{cutlocuslemma} there holds \be \lim_{\epsilon \rightarrow 0}
\int_{M-U_{\epsilon}}\Delta \psi \phi dq =\int_M \Delta \psi \phi
dq =\int_M \psi \Delta \phi dq. \ee Hence we conclude that \be
\label{laplace} \int_M  l\Delta \phi dq\leq   \liminf_{\epsilon
\rightarrow 0} \int_{M-U_{\epsilon}} \phi \Delta l dq. \ee Since
$q_0$ is arbitrary,  we deduce by using a partition of unity that
(\ref{laplaceI}) holds true for all nonnegative smooth functions
$\phi$ with compact support. \qed

\begin{lem} \label{differentialmonotonelemma}
We have on $\cup_{\tau} \Omega(\tau) \times \{\tau\}$
\ba \label{differentialequality}
l_{\tau}-\frac{R}{2}+\frac{|\nabla l|^2}{2}+\frac{l}{2\tau}=0,
\ea
 \be
\label{heatinequality} l_{\tau}-\Delta l+|\nabla l|^2
-R+\frac{n}{2\tau}\geq 0, \ee 
and \be \label{ellipticinequality}
\Delta l -\frac{|\nabla l|^2}{2}+\frac{R}{2}+\frac{l-n}{2\tau} \leq 0. \ee Moreover,
%(\ref{heatinequality}) or (\ref{ellipticinequality}) becomes an equality at a point if and only
%if \be \label{laplacel}  \Delta l=- R + \frac{n}{2 \tau}
%-\frac{1}{2\tau^{\frac{3}{2}}}K \ee holds true at that point, where $K$ is
%defined on page 16 in [P1].
(\ref{heatinequality}) becomes an equality at a point if and only 
if (\ref{ellipticinequality}) becomes an equality at that point.
\end{lem} \Pf   The equation (\ref{differentialequality}) follows from [(7.5), P1] and 
[(7.6), P1]. The inequality (\ref{heatinequality}) is [(7.13),
P1], while the inequality (\ref{ellipticinequality}) is [(7.14),
P1]. 
%The claims about the equality cases follow from the arguments
%in [P1] for [(7.13), P1] and [(7.14), P]. 
On the other hand,  the left hand side of (\ref{heatinequality}) equals the 
left hand side of (\ref{differentialequality}) minus the left hand side 
of (\ref{ellipticinequality}).  The statement about the equality cases 
follows. \qed

\begin{theo} \label{weaklemma} Assume that
the Ricci curvature is bounded from below
on $[0, \tau]$ for each $\tau \in (0, T)$. Then the equations
\be  \label{heat}
l_{\tau}-\Delta l+|\nabla l|^2-R+\frac{n}{2\tau}\geq 0
\ee
and
\be  \label{elliptic}
\Delta l -\frac{|\nabla l|^2}{2}+\frac{R}{2}+\frac{l-n}{2\tau}\leq 0
\ee
hold true on $M \times (0, T)$, when $\Delta l$ is interpreted
in the weak sense. Namely we have
\be \label{integralheat}
\int_{\tau_1}^{\tau_2} \int_M \{\nabla l \cdot
\nabla \phi+(l_{\tau}+|\nabla l|^2-R+\frac{n}{2\tau})\phi
\} dq d\tau \geq 0
\ee
for $0 < \tau_1<\tau_2<T$ and nonnegative Lipschitz functions $\phi$ on $M \times [\tau_1, \tau_2]$ with compact support, and
\be  \label{integralelliptic}
\int_M \{- \nabla l \cdot \nabla \phi + \phi(-\frac{|\nabla l|^2}{2}+\frac{R}{2}+\frac{l-n}{2\tau})\}
dq\leq 0
\ee
for nonnegative Lipschitz functions $\phi$ on $M$ with compact support and
each $\tau \in (0, T)$.  
\end{theo}
\Pf We first consider (\ref{heat}).  Let $\phi$ be a nonnegative 
Lipschitz function on $M \times [\tau_1, \tau_2]$ with compact
support. By Proposition~\ref{lipschitzlemmaI},
Proposition~\ref{lipschitzlemmaII}, Lemma \ref{measurablelemma} and
Lemma~\ref{laplacelemma} we have \ba
\label{weakinequalityargument}
%\int_0^T\int_M \{-l(\phi_{\tau}+\Delta \phi)+(|\nabla l|^2
%\nonumber\\ -R
%+\frac{n}{2\tau})\phi\} dqd\tau
%&=&
\int^{\tau_2}_{\tau_1}\int_M \{\nabla l \cdot \nabla  \phi +(l_{\tau}+|\nabla l|^2-R
+\frac{n}{2\tau})\phi\}dqd\tau
 &\geq&
\int^{\tau_2}_{\tau^1} \int_M^* (l_{\tau}-\Delta l +|\nabla l|^2-R
\nonumber\\&&+\frac{n}{2\tau})\phi dq d\tau,
\ea
where  the integral $\int_M^*$ means the limsup of the integral on
$M-U_{\epsilon}$ as $\epsilon \rightarrow 0$.
% By [(7.13),P], the inequality
%\be  \label{heatinequality}
%l_{\tau}-\Delta l+|\nabla l|^2 -R+\frac{n}{2\tau}\geq 0
%\ee
%holds true on the domain $\cup_{\tau \in (0, T)} \Omega({\tau})$. It follows
By (\ref{heatinequality}) the right hand side in
(\ref{weakinequalityargument}) is nonnegative.

The inequality (\ref{elliptic}) follows from a similar argument,
using the inequality (\ref{ellipticinequality})
%[(7.14), P]: \be \label{ellipticinequality}
%2\Delta l -|\nabla l|^2+R+\frac{l-n}{\tau} \leq 0 \ee
instead of (\ref{heatinequality}). \qed
\\

%\begin{Rk}  By the argument in the proof of Lemma \ref{concavelemma}, one easily shows
%that the differential inequalities (\ref{heat}) and (\ref{elliptic})
%hold true in the sense of barriers, hence also in the sense of viscosity, as
%is easy to verify. Using the results and arguments in [I] on general
%elliptic equations, one then obtains
%another proof of Lemma \ref{weaklemma}. One notices however that the arguments
%in [I] depend on approximations by semi-convex or semi-concave functions and
%are more involved.
%\end{Rk}

\sect{Basic Properties of the $l$-Function II}

\begin{lem} \label{minimumlemma} Assume that the Ricci curvature is bounded
from below on $[0, \tau]$.
Then the minimum of $l(\cdot, \tau)$ does not exceed
$\frac{n}{2}$.
\end{lem}
\Pf  We have the differential inequality [(7.10), P1] on $\cup\Omega(\tau) \times \{\tau\}$
\ba
\Delta L \le -2\sqrt{\tau} R+\frac{1}{\sqrt{\tau}}-\frac{1}{K}, 
\ea
which is obtained by taking trace in (\ref{hessian}). Here $K$ is
defined on page 16 in [P1]. Combining this 
with the equation [(7.5), P1] 
\ba
L_{\tau}=2\sqrt{\tau}R-\frac{1}{2\tau}L+\frac{1}{\tau}K
\ea
yields the differential
inequality [(7.15), P]: \be \label{minequation} \bar
L_{\tau}+\Delta \bar L\leq 2n,\ee 
where $\bar L(q, \tau)=2\sqrt{\tau}L(q, \tau)$. 
By the argument in the proof of Proposition
\ref{concavelemma}, one readily shows that under the assumption about the Ricci curvature
(\ref{minequation}) holds true in the sense of
barriers.  More
precisely, for each $q \in M$, each $\tau \in (0, T)$ and each
$\varepsilon>0$, there is a smooth function $u$ (an
$\varepsilon$-barrier at $(q, \tau)$) on a neighborhood of $(q,
\tau)$ in $M \times [\tau, T)$ such that $u \geq \bar L, u(q,
\tau)=\bar L(q, \tau)$ and $u_{\tau}(q, \tau)+\Delta u (q, \tau)
\leq 2n+\varepsilon.$  (We use the forward interval $[\tau, T)$
here because the left hand side of (\ref{minequation}) is the
backward heat operator.) By Lemma \ref{firstllemma} the minimum of
$l(\cdot, \tau)$ and hence of $\bar L(\cdot, \tau)$ is achieved
for every $\tau$. Consequently, the maximum principle implies that
the minimum of $\bar L(\cdot, \tau)-2n\tau$ is nonincreasing. The
desired bound for the minimum of $l$ follows.

The details of the said maximum principle are as follows. Set
$v=\bar L-2n\tau$. Then $v$ satisfies $v_{\tau}+\Delta v \leq 0$ in the
sense of barriers.  Let $h(\tau)=
\min v(\cdot, \tau)$. Consider $\tau$ and
a minimum point $q$ for $v(\cdot, \tau)$. For $\varepsilon>0$ let $u_{\varepsilon}$ be an
$\varepsilon$-barrier of $v$ at $(q, \tau)$. Then we have
for $\tau'>\tau$ sufficiently close to $\tau$
\be
\frac{h(\tau')-h(\tau)}{\tau'-\tau} \leq \frac{v(q, \tau')-v(q, \tau)}
{\tau'-\tau} \leq  \frac{u(q, \tau')-u(q, \tau)}{\tau'-\tau}.
\ee
Taking limit we obtain
\be
\frac{d^+ h}{d\tau} \le \frac{\partial u}{\partial \tau}(q, \tau)
\leq -\Delta u(q, \tau)+\varepsilon,
\ee
where $\frac{d^+ h}{d \tau}= \limsup_{\tau' \rightarrow 0^+}
\frac{h(\tau')-h(\tau)}{\tau'-\tau}$.
Obviously, $q$ is a minimum point for $u(\cdot, \tau)$, whence $\Delta u(q, \tau)
\geq 0$. Letting $\varepsilon \rightarrow 0$ we then arrive at
\be
\frac{d^+ h}{d\tau} \leq 0.
\ee
Consequently, $h$ is nonincreasing, cf. [H].
\qed
\\

Next we present a lower bound and an upper bound for
$l$ at any given time in terms of  the distance at the same time, which can be compared with Lemma \ref{firstllemma}.
The basic idea of the lower bound and its proof was communicated
to us by Perelman.  To work out the precise dependence of 
the estimate on $\tau$, we formulate it in a scaling invariant form.

\begin{lem} \label{perelmanlemma} Assume that the curvature operator
is nonnegative on $(0, T)$. Let $\bar \tau \in (0, T)$. Then we have on $(0, \bar \tau]$
\be  \label{d-l}
-l(x, \tau)-1+C_1\frac{d^2(x,q, {\tau})}{\tau} \le l(q, \tau) \le 2l(x, \tau)+C_2\frac{d^2(x, q, \tau)}{\tau}
%\le l(q, {\tau}) \le
%l(x, \tau)+C_2\frac{d^2(x,q, {\tau})}{\tau}
\ee
for all $x, q \in M$, where $C_1$ and $C_2$ are 
positive constants depending only on the dimension $n$ and the magnitude of $\frac{\bar \tau}{T-\bar \tau}$.
(In particular, $C_1$ and $C_2$ depend only on $n$ if $T=\infty$.)
\end{lem}
\Pf It follows from the Lipschitz estimate (\ref{lipschitzestimate}) in Theorem \ref{3estimateslemma} that 
\ba \label{squareestimate}
\sqrt{l}(q, \tau) \le \sqrt{l}(x, \tau) + \sqrt{\frac{C}{4\tau}}d(x, q, \tau).
\ea
Squaring it we arrive at the upper bound in (\ref{d-l}). 
Next we derive the lower bound.  Note that $l$ and the quantity $d^2(x, q, \tau)/\tau$ are
both invariant under the rescaling $g(\tau) \rightarrow a^{-1}g(a \tau)$.
Hence it suffices to prove (\ref{d-l}) for the case $\tau=1$.

Since $\Omega({1})$ is dense in $M$, it suffices to
consider the case $x, q \in \Omega({1})$.
Let $\gamma_x, \gamma_q$ be the minimal ${\cal L}_{0, 1}$-geodesics from $p$ to $x,q$
respectively.
Then
\ba
d(x, q, 1) & = & \int_0^{{1}} \frac{d}{ds} d(\gamma_x(s), \gamma_q(s), s)ds
\nonumber\\
& = &  \int_0^{{1}} \left[ \frac{\partial d}{\partial s} (\gamma_x(s),
\gamma_q(s), s) +
\nabla_I d \cdot \gamma_x' (s) + \nabla_{II} d \cdot \gamma_q' (s) \right]
ds,  \label{d-x-q}
\ea
where $\nabla_I$ refers to the gradient with respect to the first argument,
and $\nabla_{II}$ that with respect to the second argument.

By (\ref{gradientformula}) we have $\gamma_x' (s) = \nabla l(\gamma_x(s), s)$ and
$\gamma_q'(s)=\nabla l(\gamma_q(s), s)$. Since the scalar curvature is nonnegative, there holds
\be \label{now}
l(\gamma_q(s), s) = \frac{1}{2\sqrt{s}}  {\cal L}_{0, s}(\gamma_q|_{[0, s]}) \le
\frac{1}{2\sqrt{s}}  {\cal L}_{0, 1}(\gamma_q) =
{\frac{1}{\sqrt{s}}} l(q, {1}).
\ee
Similarly, we have
\ba \label{now1}
l(\gamma_x(s),s) \le {\frac{1}{\sqrt{s}}} l(x, {1}).
\ea
Hence we can apply (\ref{nablaestimate}) to deduce
\be
|\gamma_q' (s)| \le \sqrt{C}  s^{-3/4} \sqrt{l}(q, {1}),
|\gamma_x' (s)| \le \sqrt{C}  s^{-3/4} \sqrt{l}(x, {1}).
\label{gamma-q}
\ee

Next we estimate $\frac{\partial}{\partial s} d(\gamma_x(s), \gamma_q(s), s)$.
%It follows from  (\ref{nablaestimate}) that
%\be
%|\nabla l^{\frac 12}|^2 \le \frac{C}{4\tau}.  \label{gradient}
%\ee
Set $r_1(s) = s^{3/4} (l(q,{1})+1)^{-1/2}$.  By  
(\ref{squareestimate}) and (\ref{now}) we have 
for $\bar q$ with
$d(\bar q,\gamma_q(s), s) \leq r_1(s)$
\be
\sqrt{l}(\bar q, s) \le \sqrt{l}(\gamma_q(s),s) +
\frac{\sqrt{C}}{2\sqrt{s}} r_1(s) \le
(s^{-1/4}+\frac{\sqrt{C}}{2}) \sqrt{l(q, {1})+1}.
\ee
By (\ref{Restimate}) we then infer
\be \label{R1}
R(\bar q,s) \le C s^{-1}   (s^{-1/4}+\frac{\sqrt{C}}{2})^2
(l(q, {1})+1).
\ee
Similarly, we have
\be \label{R2}
R(\bar q,s) \le C s^{-1}  (s^{-1/4}+\frac{\sqrt{C}}{2})^2
(l(x, {1})+1)
\ee
for $\bar q$ with $d(\bar q,\gamma_x(s),s) \leq r_2(s)$, where
$r_2(s)=s^{3/4} (l(x,{1})+1)^{-1/2}$.

We set $r_0(s)=s^{3/4} (l(q,{1})+l(x,1)+1)^{-1/2}$. Applying [P, (8.3 (b))] to 
the present situation of the backward Ricci flow  we obtain 
\ba
 \frac{\partial }{\partial s} d(\gamma_x(s), \gamma_q(s), s) \le 
2(n-1)\left(
\frac{2}{3}Kr_0(s) + r_0(s)^{-1}\right), 
\ea
where $(n-1)K$ is an upper bound for the Ricci curvature at time 
$s$ on the geodesic balls $d(\gamma_q(s),\cdot, s) \leq r_0(s)$ and 
$d(\gamma_x(s),\cdot, s) \leq r_0(s)$. 
By the estimates 
(\ref{R1}) and (\ref{R2}) and the nonnegativity of the Ricci curvature, 
we can choose 
\ba
K=\frac{C}{n-1} s^{-1}   (s^{-1/4}+\frac{\sqrt{C}}{2})^2
(l(q, {1})+l(x, {1})+1).
\ea
Hence we deduce
\ba \label{partial-d}
\frac{\partial }{\partial s} d(\gamma_x(s), \gamma_q(s), s) &\le&
  \frac{4C}{3}\left((1+\frac{3(n-1)}{2C})s^{-3/4}+\sqrt{C} s^{-1/2}
  + \frac{C}{4}s^{-1/4} \right) \cdot \nonumber \\
&& \sqrt{l(q, {1})+l(x,1)+1}. 
\ea
Combining (\ref{d-x-q}), (\ref{gamma-q}) and
(\ref{partial-d}) we arrive at 
\ba
d(x,q,1) &\le& 4\sqrt{C}\sqrt{l}(q,1)+4\sqrt{C}\sqrt{l}(x,1)+ \nonumber \\
&&\frac{4C}{3}\left(4(1+\frac{3(n-1)}{2C})+2\sqrt{C}+\frac{C}{3}\right) 
\sqrt{l(q,1)+l(x,1)+1}  \nonumber \\
&\le& (17C^2+8(n-1))\sqrt{l(q,1)+l(x,1)+1}.
\ea
(We may assume that $C \ge 1$.)
This estimate yields
the lower bound in (\ref{d-l}) for $\tau=1$.
\qed

\begin{theo} \label{limitweaklemma}
Assume either that the curvature operator is nonnegative for each
$\tau$ or that the sectional curvature is bounded on $[0, \tau]$ for
each $\tau$. Then the inequality (\ref{integralheat}) holds true
for all $0<\tau_1<\tau_2$ and nonnegative locally Lipschitz functions $\phi$ on $M \times [\tau_1,
\tau_2]$ such that $\phi \le \bar Ce^{-\bar c l}$ and $|\nabla \phi| \le
\bar Ce^{-\bar c l}$ for positive constants $\bar C$ and $\bar c$ depending on  
$\phi$ and the magnitude of $\tau_1^{-1}$.  Similarly, the inequality (\ref{integralelliptic}) holds
true for nonnegative locally Lipschitz functions $\phi$ on $M$ such that $\phi
\le \bar Ce^{-\bar c l}$ and $|\nabla \phi| \le \bar C e^{-\bar c l}$ with positive 
constants $\bar C$ and $\bar c$ 
depending on $\phi$ and the magnitude of $\tau^{-1}$.    In both cases, the involved integrals are
absolutely convergent.  In particular, we obtain by choosing
$\phi=\tau^{-\frac{n}{2}}e^{-l}$ in (\ref{integralheat}) \be \label{monotoneintegral}
\int_{\tau_1}^{\tau_2} \int_M (l_{\tau}
-R+\frac{n}{2\tau})e^{-l}\tau^{-\frac{n}{2}} dqd\tau \geq 0. \ee
\end{theo}
\Pf  We present the case of (\ref{integralheat}), while the case
of (\ref{integralelliptic}) is similar and easier. We can assume that $M$ is
noncompact. 

 Let $0<\tau_1<\tau_2$ and $\phi$ as specified in the statements of the theorem 
 be given.  \\

 \noindent {\bf Part 1} {\it Absolute Convergence} \\
 
 We first show that the integral on the left hand side of (\ref{integralheat}) converges absolutely. 
 Indeed, we can take abolute value of every term in the integrand and still have convergence.
 In the case of bounded sectional curvature we apply Proposition \ref{boundedcurvaturelemma} to deduce for each $\tau$ 
 \ba \label{convergencebound1}
 \int_{\tau_1}^{\tau_2}\int_M (|\nabla l | \cdot  |\nabla \phi| +(|l_{\tau}|+
 |\nabla l|^2+|R|+\frac{n}{2\tau})   |\phi|) dqd\tau \le
 %\\ &\le& \int_{\tau_1}^{\tau_2} \int_M 
 %[\sqrt{\frac{C}{\tau}(l+\tau+1)} + \frac{2C}{\tau}(l+\tau+1)+|R|
 %+\frac{n}{2\tau}] \bar C e^{-\bar cl}dq d\tau \nonumber \\ &\le& 
 %\int_{\tau_1}^{\tau_2} \int_M C_1(l+1)e^{-\bar c l}dqd\tau \le 
 \int_{\tau_1}^{\tau_2} \int_M
 \tilde C e^{-\frac{\bar c}{2}l} dq  d\tau
 \ea
 for a suitable positive constant $\tilde C$ depending on $\phi$ and the magnitude of $\tau_1^{-1}$. 
 Lemma 2.3 yields 
 \ba \label{convergencebound2}
 l(q, \tau) \ge \frac{e^{-2c\tau_2}}{4\tau} d^2(p,q,\tau_1)-\frac{nc}{3}\tau_2
 \ea
 for $\tau \in [\tau_1, \tau_2]$ and a postive constant $c$ (a lower bound for the Ricci curvature). By Lemma 2.1, we have 
 $d(p,q,\tau) \le e^{C(\tau_2-\tau_1)}d(p, q, \tau_1)$ for $\tau \in [\tau_1, \tau_2]$
 with $C$ denoting an upper bound for Ricci curvature.  By this and the volume 
 comparison, we infer that $dq$ grows at most like $e^{c_1d(p,q,\tau_1)}$ for a positive 
 constant $c_1$. Hence (\ref{convergencebound1}) and (\ref{convergencebound2}) 
 yield a finite upper bound for $\int_{\tau_1}^{\tau_2} \int_M (|\nabla l| \cdot |\nabla \phi| + (|l_{\tau}|+
 |\nabla l|^2+|R|+\frac{n}{2\tau})  |\phi|) dqd\tau$.
 
 In the case of nonnegative curvature operator, we argue in a similar fashion, utilizing however 
 different lemmas. Applying Theorem \ref{3estimateslemma} we again infer 
(\ref{convergencebound1}). 
 By Lemma 2.3 or Lemma 3.2, there exists a minimum 
point $x(\tau)$ of $l(\cdot, \tau)$ for each $\tau\in(0, T)$. 
%By continuity and Lemma 3.1, there is a neighborhood 
%$I_{\tau}$ of each $\tau$ such that 
%\ba
%l(x(\tau), \tau') \le n
%\ea
%for $\tau' \in I_{\tau}$. 
By the Harnack inequality (\ref{tauharnack}) in Theorem \ref{3estimateslemma} and Lemma \ref{minimumlemma} we deduce
\ba
l(x(\tau_1), \tau) \le \frac{n}{2} (\frac{\tau_2}{\tau_1})^{C} 
\ea
for $\tau \in [\tau_1, \tau_2]$. Then Lemma \ref{perelmanlemma} leads 
to 
\ba \label{l-estimate}
l(q, \tau) \ge \frac{C_1}{\tau_2}d^2(x(\tau_1), q, \tau) -\frac{n}{2} (\frac{\tau_2}{\tau_1})^{C}-1
\ea
for $\tau \in [\tau_1, \tau_2]$ and all $q \in M$. 
%We first assume that $[\tau_1, \tau_2]$ is contained in 
%$I_{\tau^*}$ for some $\tau^*$.   By Lemma 3.2 we then deduce
%\ba \label{convergencebound3}
%l(q, \tau) \ge \frac{d^2(x(\tau^*), q, \tau)}{\tau}-n-1. 
%\ea
By volume comparison, $dq$ (at time $\tau$) grows at most at the euclidean rate, i.e. $d(x(\tau_1), q ,\tau)^{n-1}$, 
with $x(\tau_1)$ as the geodesic center.  Hence 
(\ref{convergencebound1}) and (\ref{l-estimate}) lead to a desired finite 
upper bound.  
%In the general case, we can find $\tau_1=s_1<s_2<...<s_m=\tau_2$ for 
%some $m$ along with $s_i^* \in [s_i, s_{i+1}], i=1,...,m-1$, such that 
%$[s_i, s_{i+1}]$ is contained in $I_{s_i^*}$. Applying the above argument to 
%each interval $[s_i, s_{i+1}]$ in the role of $[\tau_1, \tau_2]$ and 
%summing up, we arrive at the desired finite upper bound. 
 \\

\noindent {\bf Part 2} {\it The Integral Inequality} \\

For each natural number $k$
we choose a smooth nonnegative function $\zeta_k \le 1$ on the real line 
such that $\zeta_k=1$ on $[0, k]$, $\zeta_k=0$ on $[k+2, \infty)$ and $|\zeta_k'|\le 1$ everywhere.  In the case of bounded 
sectional curvature we define
$\eta_k$ on $M \times [\tau_1, \tau_2]$ by the formula $\eta(q, \tau)=
\zeta_k(d(p, q, \tau))$. Then we have $0\le \eta_k\le 1, |\nabla \eta_k| \le 1$ everywhere,  
$\eta_k(q,\tau)=1$ whenever $d(p, q, \tau) \le k$, and  $\eta_k(q,\tau)=0$ whenever $d(p, q, \tau) \ge k+2$.  By Theorem 
\ref{weaklemma} 
we then have 
\be \label{approximate} \int_{\tau_1}^{\tau_2} \int_M
\{\nabla l \cdot (\nabla \eta_k \phi+\eta_k \nabla
\phi)+(l_{\tau}+|\nabla l|^2-R+\frac{n}{2\tau})\eta_k\phi \} dq
d\tau \geq 0 \ee
Let
$I_k$ denote the integral $\int_{\tau_1}^{\tau_2} \int_M
\{\nabla l \cdot \eta_k \nabla
\phi+(l_{\tau}+|\nabla l|^2-R+\frac{n}{2\tau})\eta_k\phi \} dq
d\tau$, and $I$ denote the integral 
$\int_{\tau_1}^{\tau_2} \int_M
\{\nabla l \cdot \nabla
\phi+(l_{\tau}+|\nabla l|^2-R+\frac{n}{2\tau})\phi \} dq
d\tau$. Then we have
\ba 
|I_k-I| \le \int_{\tau_1}^{\tau_2} \int_{B_{k+2}(p,\tau)-B_k(p, \tau)} (|\nabla l
\cdot \nabla \phi|+|l_{\tau}+|\nabla l|^2-R+\frac{n}{2\tau}|\phi) dqd\tau.
\ea
By the above arguments for the absolute convergence in the case of bounded sectional curvature we infer 
\ba \label{Iestimate}
|I_k-I|\le \bar C_1 (\tau_2-\tau_1)e^{-\bar c_1 k}
\ea
for some positive constants $\bar C_1$ and $\bar c_1$.  Similarly, we have 
\ba \label{anotherestimate}
|\int_{\tau_1}^{\tau_2}\int_M \nabla l \cdot \nabla \eta_k \phi 
dq d\tau| \le \int_{\tau_1}^{\tau_2}\int_{M-B_k(p, \tau)} |\nabla l| \phi 
dq d\tau \le \bar C_2(\tau_2-\tau_1) e^{-\bar c_2k}
\ea
for some positive constants $\bar C_2$ and $\bar c_2$. 
%Taking limit as $k \rightarrow \infty$
%we then deduce the desired integral inequality $I \ge 0$. 

In the case of nonnegative curvature operator we 
%consider the subdivision $\tau_1=s_1<s_2<...
%s%_m=\tau_2$ and the points $s_i^*$ given before. For each fixed $i$, we define 
set $\eta_k(q, \tau)=\zeta_k(d(x(\tau_1), q, \tau))$ and use $x(\tau_1)$ as the geodesic 
center.
% on $M \times [s_i, s_{i+1}]$ (or rather on 
%$M \times [s_i, s_{i+1})$). 
Arguing as before
%above with $[s_i, s_{i+1}]$ playing the role of $[\tau_1, \tau_2]$ and 
%employing the arguments for the absolute convergence in the case of 
%nonnegative curvature operator, 
%and then summing up over $i$ 
we again obtain the estimates  
(\ref{Iestimate}) and (\ref{anotherestimate}). 

Taking limit as $k \rightarrow \infty$ we 
arrive at the desired integral inequality.
\qed \\

%\noindent {\bf Remark} In Part 1 of the proof, one can also 
%use $x(\tau)$ for each $\tau$ instead of $x(s_i^*)$.  The subdivision of $[\tau_1, \tau_2]$ 
%along with the points $x(s_i^*)$ are necessary in Part 2, because $x(\tau)$ may not 
%depend continuously  
%on $\tau$, and the function $\eta_k(d(x(\tau), q, \tau)$ (in the role of 
%$\zeta_k(d(x(s_i^*), q, \tau)$)  may not be 
%measurable.  This remark is also relevant for the proof of Theorem 5.1 below.
%\qed

\sect{The Reduced Volume}

We continue with the solution $(M, g=g(\tau))$ of the backward
Ricci flow on $[0, T)$ as before (assuming that $g(\tau)$ is
complete for each $\tau \in [0, T)$). \\

\noindent {\bf Definition 5} We define the {\it reduced volume}
(of Perelman)
$\tilde V(\tau)$ to be \be \label{volume} \tilde V(\tau)=\tilde
V_g(\tau)=\int_M \tau^{-\frac{n}{2}}e^{-l(q, \tau)} dq. \ee

A basic property of $\tilde V$ is its invariance under the rescaling 
$g(\tau) \rightarrow g_a(\tau)=a^{-1}g(a\tau)$, which easily follows from 
Lemma 2.2.  Our main goal is to obtain monotonicity of the reduced volume and its
upper bounds, and the associated rigidities. For this purpose, we need as in [P1] the following
weighted monotonicity of the Jacobian of the $\cal{L}$-exponential
map given in [P1].

\begin{lem} \label{jacobilemma} Let $J(\tau)(v)=J_g(\tau)(v)$ denote the Jacobian of the
$\cal{L}$-exponential map $\lexp$ at $v \in \Omega^{T_pM}(\tau)$,
where $T_pM$ is equipped with the metric $g(\tau)_p$. Then we have
\ba \frac{d}{d\tau} \tau^{-\frac{n}{2}}e^{-l(v, \tau)}J(\tau)(v)
\leq 0 \ea for each $v \in \Omega^{T_pM}(\tau)$, where $l(v,
\tau)= l(\gamma_v(\tau), \tau)$ ($\gamma_v$ is given in Definition
2).
 Moreover, if $\tau_1^{-\frac{n}{2}}e^{-l(v,
\tau_1)}J(\tau_1)(v) =\tau_2^{-\frac{n}{2}}e^{-l(v,
\tau_2)}J(\tau_2)(v)$ for $\tau_1 < \tau_2$ and $v \in
\Omega^{T_pM}(\tau_2)$, then the equation \ba Ric-\frac{1}{2\tau}g
+\nabla^2 l=0 \ea holds true along $\gamma_v$ on the interval
$[\tau_1, \tau_2]$.
\end{lem}
\Pf This follows from the arguments on pages 16 and 17 in [P1].
\qed
\\

\begin{lem} \label{jacobilimitlemma} Consider $v \in
\Omega(\hat \tau)$ for some $\hat \tau$. Then \be
\label{limitjacobi} \lim_{\tau \rightarrow
0}\tau^{-\frac{n}{2}}e^{-l(v, \tau)}J(\tau)(v)
=e^{-\frac{|v|^2}{4}}.  \ee Consequently, \be
\label{jacobiinequality} \tau^{-\frac{n}{2}}e^{-l(v,
\tau)}J(\tau)(v)\leq e^{-\frac{|v|^2}{4}} \ee for each $\tau$.
\end{lem}
\Pf Set $\tilde J(\tau)(v)=\tau^{-\frac{n}{2}}e^{-l(v,
\tau)}J(\tau)(v)$. The following transformation formula is easy to
verify: \be \tilde J_{g_a}(\tau)(av)dv_{g_a(0)}=\tilde J(a
\tau)(v)dv, \ee where $g_a(\tau)=a^{-1}g(a \tau)$ and $a \tau \leq
\hat \tau$.  In particular \be \label{transform} \tilde
J_{g_a}(1)(av)dv_{g_a(0)}=J(a)(v)dv. \ee  Using $exp^{{\cal{L}},
\hat \tau}$ we pull back $g, 0 \leq \tau \leq \bar \tau$ to
$\Omega^{T_p}(\bar \tau)$, and then pull it back by the scaling
map $\Phi_a(v')=av', v' \in T_pM$.  The resulting metrics will be
denoted by $g^*$.  Applying (\ref{transform}) to $g^*$ we deduce
that \be \tilde J(a)(v)=\tilde J_{g^*_a}(1)(v)\ee for $0<a<\hat
\tau$.  Next observe that over $[0, 2]$, $g^*_a$ converge smoothly
on compact sets of $T_pM$ to the euclidean steady soliton $
g^0(\tau) \equiv  g(0)_p$ as $a \rightarrow 0$. Moreover, the
image of the minimal $\cal{L}$-geodesic from the reference point
$0$ to $v$ remains in a fixed compact set during the convergence,
which follows from the arguments in the proof of Proposition
\ref{approachlemma}.  It follows that $\lim_{a \rightarrow
0}\tilde J(a)(v) = \tilde J_{g^0}(1)(v)=e^{-\frac{|v|^2}{4}}$.

The inequality (\ref{jacobiinequality}) follows from
(\ref{limitjacobi}) and Lemma \ref{jacobilemma}. \qed

\begin{theo} \label{newupper} Assume that the Ricci curvature is bounded from below
on $[0, \bar \tau]$ for some $\bar \tau$. Then $\tilde V(\tau)
\leq (4\pi)^{\frac{n}{2}}$ for each $\tau \in (0, \bar \tau]$. 
%The
%same holds if the Ricci curvature is bounded from above on $[0
%,\bar \tau]$.
\end{theo}
\Pf  By Lemma \ref{cutlocuslemma} and Lemma \ref{jacobilimitlemma}
we have \be \tilde V(\tau)= \int_{\Omega^{T_p}(\tau)}
\tau^{-\frac{n}{2}} e^{-l(v, \tau)} J(\tau)(v)dv \leq
\int_{\Omega^{T_p}(\tau)} e^{-\frac{|v|^2}{4}} dv \leq \int_{T_pM}
e^{-\frac{|v|^2}{4}} dv=(4\pi)^{\frac{n}{2}}. \ee

%An alternative argument goes as follows, for which we assume that
%the Ricci curvature is bounded on $[0, \tau]$ for each $\tau$.
%To
%obtain the desired upper bound for $\tilde V(\tau)$, we compute
%
%its limit as $\tau \rightarrow 0$.
%For $\bar \tau>0$, we set
%$g_{\bar \tau}(\tau)=\frac{1}{\bar \tau}g(\bar \tau \tau).$ By the
%scaling invariance, we have $\tilde V(\bar \tau)= \tilde
%V_{g_{\bar \tau}}(1)$.  Utilizing Lemma \ref{firstllemma},
%(\ref{basicboundI}), (\ref{basicboundII}), the volume comparison
%and the fact that $g_{\bar \tau}$ converges to the euclidean
%metric on ${\bf R}^n$ on compact sets one readily shows that the
%latter converges to $\tilde V_{g_{euc}}(1)=(4\pi)^{\frac{n}{2}}$
%as $\bar \tau \rightarrow 0$, where $g_{euc}$ denotes the
%euclidean metric on ${\bf R}^n$. It follows that $\tilde V(\tau)
%\le (4\pi)^{\frac{n}{2}}$.\qed

\begin{theo} \label{4pilemma} Assume that the Ricci curvature is nonnegative
for $s \in [0, \tau]$. Then $\tilde V(\tau) <
(4\pi)^{\frac{n}{2}}$ unless $(M, g(0))$ is isometric to ${\bf
R}^n$ and $g(s)=g(0)$ for all $s \in [0, \tau]$, in which case
$\tilde V(\tau)=(4\pi)^{\frac{n}{2}}$.
\end{theo}
\Pf   By (\ref{Ricciflow}), we have $\frac{\partial}{\partial t}
dq=Rdq$. This and (\ref{firstLestimate}) imply that \be
\label{vreduction} \tilde V(\tau) \le \tau^{-\frac{n}{2}} \int_M
e^{-\frac{d^2(p, q, 0)}{4\tau}}dq \le \tau^{-\frac{n}{2}} \int_M
e^{-\frac{d^2(p, q, 0)}{4\tau}}dq|_{0}. \ee By volume comparison,
we have \be \label{veuc} \int_M e^{-\frac{d^2(p, q,
0)}{4\tau}}dq|_{0}\le \int_{{\bf
R}^n}e^{-\frac{|x|^2}{4\tau}}dx=(4\pi \tau)^{\frac{n}{2}}. \ee
Hence we arrive at the desired inequality. If $\tilde
V(\tau)=(4\pi)^{\frac{n}{2}}$, then (\ref{veuc}) must be an
equality, and hence $(M, g(0))$ is isometric to ${\bf R}^n$. The
second inequality in (\ref{vreduction}) must also be an equality.
Consequently, $R\equiv 0$ and hence $Ric\equiv 0$ for $s\in [0,
\tau]$. It follows that $g(s)=g(0)$ for $s \in [0, \tau]$. \qed

\begin{theo} \label{monotonelemma} If the Ricci curvature
is bounded from below on $[0, \tau]$ for each $\tau$, then $\tilde
V(\tau)$ is a nonincreasing function. 
%The same holds if the Ricci
%curvature is bounded from above on $[0, \tau]$ for each $\tau$.
\end{theo}
\Pf By Lemma \ref{cutlocuslemma} we have \be \tilde
V(\tau)=\int_{\Omega(\tau)} \tau^{-\frac{n}{2}}e^{-l(q, \tau)} dq=
\int_{\Omega^{T_p}(\tau)} \tau^{-\frac{n}{2}}e^{-l(v,
\tau)}J(\tau) dv, \ee where
%$l(v, \tau)=l(\gamma_v(\tau), \tau)$,
$dv$ denotes the euclidean volume form on $T_pM$ determined by
$g(\tau)_p$.
%and $J(\tau)$ the Jacobian of $\lexp$.
By Lemma \ref{geodesiclemma} we have for $\tau_1<\tau_2$ the
inequality ${\tilde V(\tau_2)-\tilde V(\tau_1)} \le
\int_{\Omega({\tau_2})} \tau_2^{-\frac{n}{2}}e^{-l(v,
\tau_2)}J(\tau_2) dv -\int_{\Omega(\tau_2)}
\tau_1^{-\frac{n}{2}}e^{-l(v, \tau_1)}J(\tau_1) dv$. By Lemma
\ref{jacobilemma} we then obtain the desired monotonicity. \qed

\begin{lem} \label{volumedifferenceprop}
 Assume
either that the sectional curvature is bounded on $[0, \tau]$ for each
$\tau$ or that the curvature operator is nonnegative. Then there
holds \ba \label{volumedifference} {\tilde V(\tau_2)-\tilde
V(\tau_1)}&=&-\int_{\tau_1}^{\tau_2}\int_M (l_{\tau}
-R+\frac{n}{2\tau})e^{-l}\tau^{-\frac{n}{2}} dqd\tau 
%&=&\int_{\tau_1}^{\tau_2} \int_M R
%e^{-l}\tau^{-\frac{n}{2}}dq d\tau +\int_M e^{-l(q,
%\tau_2)}dq|_{\tau_2} \nonumber\\ &&-\int_M e^{-l(q,
%\tau_1)}dq|_{\tau_1}
\ea
for all
$0<\tau_1<\tau_2 <T$.
\end{lem}
\Pf We first assume nonnegative curvature operator.  Consider
$\tau_2>\tau_1$. Applying Proposition \ref{lipschitzlemmaII}
%Lemma \ref{3estimateslemma} and Lemma \ref{perelmanlemma}
 we
deduce
$$
{\tilde V(\tau_2)-\tilde V(\tau_1)}=\int_M \int_{\tau_1}^{\tau_2}
\frac{\partial}{\partial \tau} (\tau^{-\frac{n}{2}}e^{-l(q, \tau)}
dq) =
$$
\be -\int_M
\int_{\tau_1}^{\tau_2}(l_{\tau}-R+\frac{n}{2\tau})e^{-l}\tau^{-\frac{n}{2}}
dqd\tau. \ee   By the proof of Theorem \ref{limitweaklemma}, the last integral is absolutely convergent.
Hence  we can
switch the integration order to arrive at the first equation in
(\ref{volumedifference}).
% Next we have 
%\ba
%-\int_{\tau_1}^{\tau_2} \int_M (l_{\tau}+\frac{n}{2\tau})e^{-1}\tau^{-\frac{n}{2}}dqd\tau 
%=-\int_M \int_{\tau_1}^{\tau_2}  (l_{\tau}+\frac{n}{2\tau})e^{-1}\tau^{-\frac{n}{2}} d\tau dq 
%\ea

The proof of  (\ref{volumedifference}) in the case of bounded
sectional curvature is similar. Note that in this case we also have \be
\frac{\partial \tilde V}{\partial \tau}=\int_M (l_{\tau}
-R+\frac{n}{2\tau})e^{-l}\tau^{-\frac{n}{2}} \ee for every $\tau$.
This can be seen by computing the relevant difference quotient and
applying the dominated convergence theorem to pass to limit.

We remark that (\ref{volumedifference}) and Theorem
\ref{limitweaklemma} also imply the monotonicity of $\tilde
V(\tau)$. \qed
\\

\noindent {\bf Definition 6} In this definition, let $g$ be a
smooth solution of the backward Ricci flow on $N \times I$ for a
smooth manifold $N$ and an interval $I$. We say that $g$ is a {\it
gradient shrinking soliton} with {\it time origin} $\tau_0$ and
{\it potential function} $f$ on an open subset $O$ of $N \times
I$, where $f$ is a smooth function on $O$, provided that $g$
satisfies the gradient shrinking soliton equation \ba
\label{solitonequation}
Ric-\frac{1}{2(\tau-\tau_0)} g +\nabla^2 f=0 \ea in $O$.\\

\begin{lem} \label{solitonrule} Let  $g$ be a gradient shrinking soliton on $N \times I$
 with time origin $\tau_0$ and potential function $f$.
 Then $g$
evolves  by the pullback of a family of diffeomorphisms coupled
with scaling. More precisely, we have \be \label{newtauequation}
g(\tau)=\frac{\tau-\tau_0}{\bar \tau-\tau_0}(\phi^{-1})^*g(\bar \tau),
\ee where $\bar \tau $ is an arbitary point in $I$ and $\phi$ is the solution of the
equation $\frac{\partial \phi}{\partial \tau}=\nabla f(\phi)$ with
$\phi(\bar \tau)=id$ ($id$ denotes the identity map of $N$).
\end{lem}
\Pf  We have \be \frac{\partial \bar g}{\partial \tau}=2Ric_{\bar
g}=\frac{1}{\tau-\tau_0}\bar g-2\nabla_{\bar g}^2
f=\frac{1}{\tau-\tau_0}\bar g-L_{\nabla_{\bar g} f}\bar g. \ee
Hence \be \frac{\partial }{\partial \tau} \phi^* \bar
g=\frac{1}{\tau-\tau_0} \phi^* \bar g. \ee The equation
(\ref{newtauequation}) follows. \qed

\begin{lem} \label{solitontest} As in Definition 6, let $g$ be a
smooth solution of the backward Ricci flow on $N \times I$. Let
$f$ be a smooth function on an open subset $O$ of $N \times I$. We
set $u=(4\pi\tau)^{-\frac{n}{2}} e^{-f}$ and $v=[\tau(2 \nabla f -
|\nabla f|^2 +R)+f-n]$. Then we have \be \label{boxformula} \Box
v= -2\tau|Ric+\nabla^2 f -\frac{1}{\tau} g|^2u+ 2\tau u\Delta
(u^{-1} \Box u), \ee where $\Box u=u_{\tau}-\nabla u +Ru$.
Consequently, if $u$ satisfies the heat equation $\Box u=0$, or
equivalently \be \label{testequation} \frac{\partial f}{\partial
\tau}-\Delta f+ |\nabla f|^2-R+\frac{n}{2\tau}=0, \ee
 then $v$
satisfies the equation \be \Box v=-2\tau|Ric+\nabla^2 f
-\frac{1}{\tau} g|^2. \ee  In particular, if (\ref{testequation})
holds, then $g$ is a gradient shrinking soliton on $O$ with time
origin $0$ and potential function $f$ if and only if~ $\Box v=0$
in $O$.
\end{lem}
\Pf This is a reformulation of [Proposition 9.1, P]. The formula
(\ref{boxformula}) follows from routine computations, see
e.g. [KL]. \qed \\

Now we return to our previous $g$ on $M \times [0, T)$.

\begin{theo} \label{solitontheorem} Assume either that the sectional curvature is bounded on
$[0, \tau]$ for each $\tau<T$, or that the curvature operator is
nonnegative. Assume that $\tilde V(\tau_1)=\tilde V(\tau_2)$ for
some $\tau_1<\tau_2$. Then $l$ is smooth on $M \times (\tau_1,
\tau_2)$ and $g$ is a gradient shrinking soliton on $M \times
(\tau_1, \tau_2)$ with time origin $0$ and potential function $l$.
\end{theo}
\Pf Assume $\tilde V(\tau_1)=\tilde V(\tau_2)$ for some
$\tau_1<\tau_2$. By Lemma \ref{volumedifferenceprop} we have
\ba \label{integralequality} \int_{\tau_1}^{\tau_2}\int_M
(l_{\tau} -R+\frac{n}{2\tau})e^{-l}\tau^{-\frac{n}{2}} dqd\tau=0.
\ea We set \ba Q_{\tau_1, \tau_2}(\phi)=\int_{\tau_1}^{\tau_2}
\int_M \{\nabla l \cdot \nabla \phi+(l_{\tau}+|\nabla
l|^2-R+\frac{n}{2\tau})\phi \} dq d\tau \ea for admissible $\phi$,
which are locally Lipschitz functions $\phi$ on $M \times [\tau_1,
\tau_2]$ such that $|\phi| \leq C\tau^{-\frac{n}{2}} e^{-l}$ and $|\nabla \phi| \leq
C\tau^{-\frac{n+1}{2}}\sqrt{l+\tau+1}e^{-l}$ for some bound factor $C>0$ depending on $\phi$. 
By Proposition \ref{lipschitzlemmaIII}, Proposition \ref{boundedcurvaturelemma} and Theorem  
\ref{3estimateslemma}, the function $\tau^{-\frac{n}{2}}e^{-l}$ is admissible. 
By (\ref{integralequality}) we have  $Q_{\tau_1,
\tau_2}(\tau^{-\frac{n}{2}}e^{-l})=0$. For an arbitrary
nonnegative admissible $\phi$ with bound factor $C$  we have by
Theorem \ref{limitweaklemma} $Q_{\tau_1, \tau_2}(\phi) \ge 0$ and $Q_{\tau_1,
\tau_2}(C\tau^{-\frac{n}{2}}e^{-l}-\phi) \geq 0$, whence
$0 \le Q_{\tau_1, \tau_2}(\phi) \leq
CQ_{\tau_1, \tau_2}(\tau^{-\frac{n}{2}}e^{-l})=0$, i.e. $Q_{\tau_1,
\tau_2}(\phi)=0$. By linearity of $Q_{\tau_1, \tau_2}$ we then
infer that $Q_{\tau_1, \tau_2}(\phi)=0$ for all admissible $\phi$
(simply write $\phi$ as the sum of its positive and negative parts),
in particular for all Lipschitz $\phi$ with compact support.  The
standard regularity theory for parabolic equations implies that
$l$ is smooth on $M \times (\tau_1, \tau_2)$ and satisfies \be
\label{limitheat} \frac{\partial l}{\partial \tau}-\Delta l+
|\nabla l|^2-R+\frac{n}{2\tau}=0. \ee By Lemma
\ref{differentialmonotonelemma} we then also have \be
\label{limitelliptic} 2\Delta l -|\nabla l|^2+ R+\frac{l-n}{\tau}
=0. \ee

Now we can apply Lemma \ref{solitontest} with $f=l$. By
(\ref{limitheat}), the equation (\ref{testequation}) holds true.
By (\ref{limitelliptic}), $v=0$. Hence we conclude that $g$ is a
gradient shrinking soliton with time origin $0$ and potential
function $l$.
 (The implication of Lemma
\ref{solitontest}, i.e. [9.1, P]  was first pointed out to us by
G.~Wei. Note that a similar argument is used in the proof of
Theorem 10.1 in [P1].)
\qed \\

%\noindent {\bf Remark } The converse of this result is also true.
%Using a different argument we can show that the same conclusions
%hold under the assumption that the Ricci curvature is bounded from
%below or bounded from above. (Hence Theorem \ref{upperbound2}
%below also holds under the assumption of lower bounds for the
%Ricci curvature.) These will be presented elsewhere. Although the
%argument here does not yields the most general result in the
%context of Theorem \ref{solitontheorem}, its independent interest
%merits a detailed presentation. A notable application of this
%argument can be found in the proof of Theorem
%\ref{theorem1} below. \\

\begin{theo} \label{upperbound2} Assume that the sectional curvature
is bounded on $[0, \tau]$. Then $\tilde V(\tau) <
(4\pi)^{\frac{n}{2}}$ unless $(M, g(0))$ is isometric to ${\bf
R}^n$ and $g(s)=g(0)$ for each $s \in [0, \tau]$.
\end{theo}
\Pf By Theorem \ref{newupper}, $\tilde V(\tau) \leq
(4\pi)^{\frac{n}{2}}.$ Assume that the equality holds. By Theorem
\ref{monotonelemma} and Theorem \ref{solitontheorem}, $g$ is a
gradient shrinking soliton on $M \times (0, \tau)$ with time
origin $0$ and potential function $l$. By Lemma \ref{solitonrule},
$g(\tau')= \frac{\tau'}{\tilde \tau} \phi^* g(\tilde \tau)$ for $
\tau', \tilde \tau \in (0, \tau)$. Since the sectional curvature is bounded,
we can let $\tilde \tau \rightarrow 0$ to deduce that $g(\tau')$ is
flat for each $\tau' \in (0, \tau)$. The desired
conclusion then follows from Theorem \ref{4pilemma}. \qed

%Next assume that $\tilde V(\tau_1)=\tilde V(\tau_2)$. Then the
%above arguments imply that $\Omega(\tau)=\Omega(\tau_1)$ is for
%all $\tau \in [\tau_1, \tau_2]$, and that $g$ is a gradient
%shrinking soliton with potential function $l$ on $\Omega(\tau_1)
%\times [\tau_1, \tau_2]$.

 %The conclusion in the equality case also
%follows from Lemma \ref{jacobianlemma}.
%[(7.12), P], $\tau_2^{-\frac{n}{2}}e^{-l(v,
%\tau_2)}J(\tau_2) \le \tau_1^{-\frac{n}{2}}e^{-l(v,
%\tau_1)}J(\tau_1)$. It follows that $\tilde V(\tau_2)\le \tilde
%V(\tau_1)$.

%Next we assume that $g$ is an integrable gradient shrinking
%soliton with potential function $l$ on $(0, \bar \tau]$. Then

%\begin{theo}
%Assume that the Ricci curvature is bounded on $[0, \tau]$ for each
%$\tau$. Then $\tilde V(\tau) \leq (4\pi)^{\frac{n}{2}}$ for each
%$\tau$. The equality $\tilde V(\tau)=(4\pi)^{\frac{n}{2}}$ (for
%any given $\tau$) holds true if and only if $l$ is smooth and  $g$
%is a gradient shrinking soliton on $M \times (0, \tau]$ with
%potential function $l$.
%\end{theo}

Department of Mathematics,
 University of California,
Santa Barbara, CA 93106

yer@math.ucsb.edu

\end{document}